\documentclass[11pt]{article}
\usepackage[usenames,dvipsnames]{color}
\usepackage{amsmath,amsfonts, amssymb,graphicx,geometry,mathtools,esint,enumerate,color,soul,ulem,caption,subcaption} 
\usepackage [english]{babel}
\usepackage [autostyle, english = american]{csquotes}
\MakeOuterQuote{"}

\DeclareMathOperator{\Tr}{Tr}
\DeclareMathOperator{\adj}{adj}
\DeclareMathOperator{\cof}{cof}

\DeclareMathOperator{\id}{id}
\DeclareMathOperator{\rank}{rank}

\DeclareMathOperator{\dist}{dist}
\DeclareMathOperator{\divg}{div}
\DeclareMathOperator{\curl}{curl}

\DeclareMathOperator{\range}{range}

\DeclareMathOperator{\supp}{supp}

\newcommand{\calK}{{\mathcal K}}

\newtheorem{theorem}{Theorem}[section]

\newtheorem{lemma}[theorem]{Lemma}
\newtheorem{proposition}[theorem]{Proposition}

\newtheorem{remark}[theorem]{Remark}

 \def\endproof{\hspace*{\fill}\mbox{\ \rule{.1in}{.1in}}\medskip }
\definecolor{bblue}{rgb}{0.0,0.0,.65}

\numberwithin{equation}{section}

\title{Effective behavior of nematic elastomer membranes}
\date{}
\author{
Pierluigi Cesana\\
{\small Mathematical Institute,  Woodstock Road, Oxford OX26GG, England}\\  \\
Paul Plucinsky and Kaushik Bhattacharya\\
\small California Institute of Technology, Pasadena, CA 91125, USA}

\begin{document}
\maketitle

\begin{center}
\textit{Dedicated to Jerald L. Ericksen on the occasion of his $90^{\textrm{th}}$ birthday}\\
\textbf{In press on the Archive for Rational Mechanics and Analysis}
\end{center}

\begin{abstract}
We derive the effective energy density of thin membranes of liquid crystal elastomers as the $\Gamma-$ limit of a widely used bulk model.  These membranes can display fine-scale features both due to wrinkling that one expects in thin elastic membranes and due to oscillations in the nematic director that one expects in liquid crystal elastomers.  We provide an explicit characterization of the effective energy density of membranes and the effective state of stress as a function of the planar deformation gradient.  We also provide a characterization of the fine-scale features.   We show the existence of four regimes: one where wrinkling and microstructure reduces the effective membrane energy and stress to zero, a second where wrinkling leads to uniaxial tension, a third where nematic oscillations lead to equi-biaxial tension and a fourth with no fine scale features and biaxial tension.  Importantly, we find a region where one has shear strain but no shear stress and all the fine-scale features are in-plane with no wrinkling.
\end{abstract}

\tableofcontents

\section{Introduction}

Liquid Crystal Elastomers are rubber-like solids that display unusual mechanical properties like soft elasticity and develop fine-scale microstructure under deformation.  This material consists of cross-linked polymer chains where rigid rod-like elements (mesogens) are either incorporated into the main chain or are pendent from them.   These mesogens have temperature-dependent interaction which results in phases of orientational and positional order \cite{DP93, WT03}. We refer to two phases: a high temperature {\it isotropic} phase, where thermal fluctuations thwart any attempt at order and a {\it nematic} phase, where the mesogens have a characteristic orientation but no positional order.  This average orientation of the mesogens in the nematic phase is represented by a director.

Nematic-elastic coupling is a key feature of these materials \cite{WBT94, O94}.  The isotropic to nematic phase transformation is accompanied by a very significant distortion of the solid: typically elongation along the director and contraction transverse to it.  Further, the director can rotate relative to the polymer matrix.  This novel mechanism induces a degeneracy in the low energy states associated with the entropic elasticity of the polymer network, whereby the material has a non-trivial set of nearly stress-free shape changing configurations.   This degeneracy can lead to fine-scale microstructure like stripe domains where the director alternates between two orientations in alternating stripes.  Together, all of this gives rise to soft-elasticity \cite{WT03}.

A theory of nematic elastomers, and specifically the entropic elasticity associated with it, was formulated in Warner et al. \cite{WBT94}, and was used to show the emergence of stripe domains and soft-elasticity.  Mathematically, the energy functional is not weakly lower-semicontinuous resulting in possible non-existence of minimizers: briefly minimizing sequences develop rapid oscillations that result in lower energy than its  weak limit.  These rapid oscillations are interpreted as the fine-scale microstructure in the material.  DeSimone and Dolzmann \cite{DD00}  computed the relaxation wherein the energy density is replaced with an effective energy density that accounts for all possible microstructures.  The effective energy does indeed show soft elasticity, and can be used  as by Conti et al. \cite{CDD02} to explain complex deformation patterns in clamped stretch experiments on nematic elastomer sheets \cite{KF95}.

Experiments on nematic elastomers, like the one highlighted, have largely been performed on thin sheets or membranes. These structures typically have instabilities such as wrinkling, and consequently membranes of usual elastic materials are unable to sustain compression and the state of stress is limited to uniaxial and biaxial tension.  Thus elastic membranes have been described heuristically by theories like the tension field theory of Mansfield \cite{M70},  and such theories have been obtained systematically from three-dimensional theories \cite{Pip86,SP89,LR95}.

The goal of this work is to derive an effective theory of thin membranes of liquid crystal elastomers that accounts not only for the formation of fine-scale microstructure but also instabilities like wrinkling.    An important insight that results from this is the possible states of stress in these materials.  We find that like usual elastic membranes, membranes of liquid crystal elastomers are also incapable of sustaining compression, and the state of stress is limited to uniaxial and biaxial tension.  Importantly, due to the ability of these materials to form microstructure, there is a large range of deformation gradients involving unequal stretch where the state of stress is purely equi-biaxial.   Consequently, a membrane of this material has zero shear stress even when subjected to a shear deformation within a certain range.

We start with a three dimensional variational model of liquid crystal elastomers, derive the effective behavior of a membrane -- a domain where one dimension is small compared to the other two -- as the $\Gamma-$limit of a suitably normalized functional as the ratio of these dimension goes to zero following LeDret and Raoult \cite{LR95} and others \cite{BJ99,Shu00,CD06}.   Our variational model is based on a Helmholtz free energy density that has two contributions.  The first contribution captures the elasticity associated with the polymer matrix.  Developed by Bladon, Terentjev and Warner \cite{WT03,BWT93}, it is a generalization of the classical neo-Hookean model to account for the local anisotropy due to the director.  The second contribution following Frank \cite{Frank58} penalizes the spatial non-uniformity of directors, and has been widely used in the study of liquid crystals.  In the context of liquid crystal elastomers, this penalizes domain walls -- narrow regions that separate domains of uniform director.  The competition between entropic elasticity and Frank elasticity, precisely the square-root of the ratio of the moduli $\kappa$ of the Frank elasticity to $\mu$ of the entropic elasticity -- introduces a length-scale.   It turns out (e.g. \cite{WT03}) that $\sqrt{\kappa / \mu} \sim 10-100\; nm$.  Note that the thickness $h$ of a realistic membrane is on the order of $1-100\; \mu m$ depending on the application.  Thus, one has two small parameters, and one needs to study the joint limit as both $\sqrt{\kappa/\mu}$ and $h$ go to zero, but at possibly different rates.  We do so by setting $\kappa = \kappa_h$ and studying the limit $h \rightarrow 0$.   

We find in Theorem \ref{1407251605} that the $\Gamma-$limit and thus the resulting theory is independent of the ratio $\kappa_h/h$.   This is similar of the result of Shu \cite{Shu00} in the context of membranes of materials undergoing martensitic phase transitions.  In other words, the length-scale on which the material can form microstructure does not affect the membrane limit as long as it is small compared to the lateral extent of the membrane.  Consequently, the $\Gamma-$limit we obtain coincides with the result of Conti and Dolzmann \cite{CD06} who studied the case $\kappa = 0$.  In fact, our proof draws extensively from their work.  Specifically, their result provides a lower bound and our recovery sequence is adapted from theirs.

The $\Gamma-$limit is characterized by an energy per unit area that depends only on the tangential gradient of the deformation.  It is obtained from the density of the entropic elasticity by minimizing out the normal component followed by relaxation or quasi-convexification.  We compute this by obtaining upper and lower bounds, and provide an explicit formula  in Theorem \ref{1407231520} (also shown schematically in Figure 1).  It is characterized by four regions depending on the in-plane stretch: ${\mathcal S}$ a solid region where there is no relaxation, ${\mathcal L}$ a liquid region where wrinkling and microstructure formation drive the effective energy to zero, ${\mathcal W}$ a wrinkling region where wrinkling relaxes the energy and ${\mathcal M}$ a microstructure region where stripe domains relaxe the energy.    The techniques employed here are in the same spirit as those employed by DeSimone and Dolzmann \cite{DD00} in three dimensional nematic elastomers.

We also study the oscillations related to the relaxation by characterizing the gradient Young measures associated with the minimizing sequences in Theorem \ref{th:char}.  We  show that the oscillations in the region $\mathcal{M}$ are necessarily planar oscillations of the nematic director and involve no out of plane deformation while those in the region $\mathcal{W}$ are characterized by uniform nematic director and wrinkling.

We use the characterization of the gradient Young measure to define the effective state of stress, and show that this coincides with the derivative of the effective or relaxed energy in Theorem \ref{StressTheorem}.  The Cauchy stress is given in (\ref{eq:memStress}): it is general biaxial tension in $\mathcal S$, zero in $\mathcal L$, uniaxial tension in $\mathcal W$ and equi-biaxial tension in $\mathcal M$.  As described above, the unique attributes of liquid crystal elastomers give rise to this region of equi-biaxial tension compared to membranes of usual elastic materials.

This paper is organized in the following manner: In Section \ref{1407231507}, we fix some notation and comment on background results which are used throughout the paper.  In Section \ref{ModNem}, we describe our model for nematic elastomers, a model which incorporates the entropic elasticity of the polymer matrix and an elastic penalty on the spatial gradient of the director.  In Section \ref{MemTheory}, we derive our effective theory for nematic elastomer membranes based on a notion of $\Gamma$-convergence.  In Section \ref{EffEn}, we provide an explicit formula for the energy density in our effective theory.  In Section \ref{Micro}, we characterize the microstructure in the aforementioned regions $\mathcal{M}$ and $\mathcal{W}$.  Finally, in Section \ref{stress1}, we conclude with a notion of stress in this effective theory and its physical implications.

\section{Preliminaries}\label{1407231507}
 
We gather here the notation and some background results which we use throughout the paper.
We denote with $\mathbb{R}^n$ the $n$ dimensional Euclidian space endowed with the 
usual scalar product $u \cdot v = u^T v$ and norm $|u|=\sqrt{u^T u}$. The unit sphere in $\mathbb{R}^n$ is denoted by $\mathbb{S}^{n-1}$ and it is defined as the set of all vectors $u\in\mathbb{R}^n$ with $|u|=1$. 
The space of $m\times n$ matrices with real entries is labeled with 
$\mathbb{R}^{m\times n}$. When $m>1$ we denote with $O(n)$ the orthogonal group of the matrices $F\in \mathbb{R}^{n\times n}$ for which $F F^T= F^T  F=I$, where $I$ is the identity in $\mathbb{R}^{n\times n}$
and with $SO(n)$ the rotation group of the matrices $F\in O(n)$ with $\det F=1$. 
Letting now $F\in\mathbb{R}^{m\times n}$,
$\adj_s(F)$ stands for the matrix of all $s\times s$ minors of $F$, $2\leq s\leq \min\{m,n\}$.
In the case $m=3,n=2$, if
\begin{displaymath}
F=\left(
\begin{array}{ccc}
f_{11} & f_{12}   \\
f_{21} & f_{22}   \\
f_{31} & f_{32}  
\end{array} \right), \,\textrm{then}\,\,\,\adj_2(F):=
\left(
\begin{array}{ccc}
(f_{21}f_{32}-f_{22}f_{31}) \\
-(f_{11}f_{32}-f_{12}f_{31})  \\
(f_{11}f_{22}-f_{12}f_{21})
\end{array} \right).
\end{displaymath}
If $u:\mathbb{R}^2\to\mathbb{R}^3$ is a smooth map, then $\adj_2\nabla u$ is normal to the surface with equation $\{u(x): x\in\mathbb{R}^2\}$.
Letting $F\in\mathbb{R}^{3\times 2}$, $Q\in SO(3)$, $R\in O(2)$ by Proposition 5.66 \cite{Dacorogna}
it follows that
\begin{eqnarray}\label{1303161951}
|\adj_2(F)|
=|\adj_2(QFR)|.
\end{eqnarray}
Later in this paper we label the norm of $\adj_2(F)$, with $F$ a $3\times 2$ matrix, with $\delta=\delta(F):=|\adj_2(F)|$. Furtheremore, we simply write $\adj\equiv\adj_2$ both when dealing with the adjugate of $3\times 3$ and $3\times 2$ matrices.

Finally, we state a version of the polar decomposition theorem: given any $F \in {\mathbb R}^{3\times2}$ and any rectangular Cartesian basis, there exist $\lambda_1 \ge \lambda_2 \geq 0, Q \in SO(3), R \in O(2)$ such that 
\begin{align} \label{eq:pd}
F = QDR
\end{align}
for 
\begin{align} \label{eq:d} 
D = \left( \begin{array}{cc} \lambda_1 & 0 \\ 0& \lambda_2 \\ 0&0 \end{array} \right).  
\end{align}
In fact, $\lambda_1, \lambda_2$ are the principle values of $F$.
\bigskip

We now recall some concepts in calculus of variations (cf. \cite{Dacorogna}).
We say that $f:\mathbb{R}^{m\times n}\to\mathbb{R}\cup\{+\infty\}$ is polyconvex if there exists a convex function $g$ which depends on the vector $M(F)$ of all the minors of $F$ such that $f(F)=g(M(F))$. In the case $m=3, n=2$ then $f(F)=g(F,\adj(F))$ with $g:\mathbb{R}^{9}\to\mathbb{R}\cup\{+\infty\}$.
We say that $f:\mathbb{R}^{m\times n}\to\mathbb{R}\cup\{+\infty\}$ is quasiconvex if, at every $F\in\mathbb{R}^{m\times n}$, we have
$$
\int_{(0,1)^n}f(F)dx\leq \int_{(0,1)^n}f(F+\nabla u)dx
$$
for every $u\in W^{1,\infty}_0((0,1)^n,\mathbb{R}^m)$. Note that the foregoing inequality holds for every $D$ open and bounded subset of $\mathbb{R}^n$ with $|\partial D|=0$ {\cite{BallMurat}}.
Finally, $f:\mathbb{R}^{m\times n}\to\mathbb{R}\cup\{+\infty\}$ is rank-one convex if $t\to f(F+t R)$ is a convex function for all $F,R\in\mathbb{R}^{m\times n}$ with $\rank(R)=1$.

If a function $f:\mathbb{R}^{m\times n}\to\mathbb{R}\cup\{+\infty\}$ is not quasiconvex, we define $f^{qc}$ the quasiconvex envelope of $f$ as
$$
f^{qc}:=\sup\{h\leq f, h \,\,\textrm{quasiconvex}\}.
$$
Analogously, we define $f^{c}, f^{pc}, f^{rc}$ as the convex, polyconvex and rank-one convex envelopes respectively of $f$.
In the general case of extended-value functions, convexity implies polyconvexity and polyconvexity implies both rank-one convexity and quasiconvexity, but quasiconvexity alone  does not imply rank-one convexity. Therefore, if $f:\mathbb{R}^{m\times n}\to\mathbb{R}\cup\{+\infty\}$, we have
\begin{eqnarray}\label{1303082357}
f^{pc}\leq f^{qc},\quad f^{pc}\leq f^{rc}.
\end{eqnarray}
On the other hand, in the case of a real-valued functions, quasiconvexity implies rank-one convexity and hence, if $f:\mathbb{R}^{m\times n}\to\mathbb{R}$, we have
\begin{eqnarray}\label{1303082354}
f^{c}\leq f^{pc}\leq f^{qc}\leq f^{rc}.
\end{eqnarray}
We give an alternative representation formula for the rank-one convex envelope of a function $f:\mathbb{R}^{m\times n}\to\mathbb{R}\cup\{+\infty\}$
$$
f^{rc}(F):=\inf\Bigl\{\sum_i^K\lambda_i f(F_i):\sum_i^K\lambda_i F_i=F, (\lambda_i,F_i) \,\,\textrm{satisfy} \, H_K\Bigr\} 
$$
with $\lambda_i\geq 0$ and $\sum_i^K\lambda_i=1$.
Family $(\lambda_i,F_i)$ satisfies a compatibility condition here labelled with $H_K$ and defined in \cite[Sec. 5.2.5]{Dacorogna}.
In the same spirit we define semiconvex hulls of a compact set $\calK\subset\mathbb{R}^{m\times n}$. The set
$$
\calK^{pc}=\Bigl\{F\in\mathbb{R}^{m\times n}: f(F) \leq\sup_{X\in \calK}f(X)\,\textrm{ for all}\,\, f:\mathbb{R}^{m\times n}\to\mathbb{R} \,\,\textrm{polyconvex}  \Bigr\}
$$
is the polyconvex hull of $\calK$. The quasiconvex hull $\calK^{qc}$ and the rank-one
convex hull $\calK^{rc}$ are defined analogously. The lamination convex hull $\calK^{lc}$ of $\calK$ is defined
$$
\calK^{lc}=\Bigl\{F\in\mathbb{R}^{m\times n}: f(F) \leq\sup_{X\in \calK}f(X)\,\textrm{ for all}\,\, f:\mathbb{R}^{m\times n}\to\mathbb{R}\cup\{+\infty\} \,\,\textrm{rank\hbox{-}one convex} \Bigr\}.
$$
Equivalently, $\calK^{lc}$ can be defined by succesively adding rank-one segments (see \cite{DD00}), i.e.
$$
\calK^{lc}=\bigcup_{i=0}^{\infty}\calK^{(i)}
$$
where $\calK^{0}=\calK$ and
$$
\calK^{(i+1)}=\calK^{(i)}\cup\bigl\{F=\lambda F_1+(1-\lambda)F_2:F_1,F_2\in \calK^{(i)},\textrm{rank}(F_1-F_2)\leq 1,\lambda\in[0,1]\bigr\}.
$$
The relations between the different notions of convexity imply the inclusions (see \cite{DD00})
$$
\calK^{lc}\subseteq \calK^{rc}\subseteq \calK^{qc}\subseteq \calK^{pc}.
$$
We refer the interested reader to \cite{Dacorogna} and \cite{M99}
 for a discussion of all the different notions of convexity and their relations.
\bigskip

Finally, we introduce the notion of a  gradient Young measure that  characterize the statistics of the fine-scale oscillations in the gradients  weakly converging sequences (cf. \cite{M99}).  We define a homogenous $H^1$ gradient Young measure to be a probability measure that satisfies Jensen's inequality for every quasiconvex function $f:\mathbb{R}^{3\times 2}\to \mathbb{R}$ whose norm can be bounded by a quadratic function.  Let $M$ denote the space of signed Radon measures on ${\mathbb R}^{3 \times 2}$ with the finite mass paring 
\begin{align*}
\langle \mu , f \rangle = \int_{\mathbb{R}^{3\times2}} f(\tilde{G}) d \mu(\tilde{G}).
\end{align*}
Then the space of homogenous $H^1$ gradient Young measures is given by 
\begin{align} \label{eq:mqc}
M^{qc} := &\big\{ \nu  \in M :   ||\nu|| = 1,  \langle \nu , f \rangle \geq f( \langle \nu , \id \rangle ) \nonumber \\
 & \quad \quad \forall \;\; f:\mathbb{R}^{3\times2} \rightarrow \mathbb{R} \text{ quasiconvex with } |f(\tilde{G})| \leq C(|\tilde{G}|^2 +1) \big\}.  
\end{align}

\section{Model of nematic elastomers}\label{ModNem}
Consider a nematic elastomer  occupying a region $\Omega$ in its reference configuration, and assume that it is in its stress-free isotropic state in this configuration.  Let $y:\Omega \rightarrow {\mathbb R}^3$ describe the deformation and $n:\Omega \rightarrow \mathbb{S}^2$ describe the director field.   We denote $\nabla u$ to be the reference gradient of some field $u:\Omega\rightarrow {\mathbb R}^3$ and $\nabla_y u$ to be the spatial gradient of $u$.  It follows $\nabla_y u = \nabla u (F)^{-1}$ where $F= \nabla y$ is the deformation gradient.

We take the Helmholtz free energy density of the nematic elastomer to be the sum of two contributions:
$$
W = W^{e} + W^{n}
$$
where the $W^{e}$ describes the entropic elasticity of the underlying polymer chains of the nematic elastomer and $W^{n}$ describes the elasticity of the nematic mesogens.

Following Bladon, Terentjev and Warner in \cite{BWT93,WT03}, we take the entropic elasticity to be of the form
\begin{eqnarray}\label{eq:We1}
W^{e}(F,n)= \begin{cases} 
\displaystyle{\frac{\mu}{2}\Bigl(\Tr (F^T \ell^{-1} F)-3\Bigr)} &  \textrm{if } \det F=1, |n| = 1, \\
+\infty & \textrm{else,} 
\end{cases}
\end{eqnarray}
where 
\begin{equation}\label{eq:stepLength}
\ell = r^{-1/3} (I + (r-1) n \otimes n)
\end{equation}
is the step-length tensor.  
Here $\mu>0$ is the shear modulus of rubber and $r \geq 1$ is the (non-dimensional) backbone anisotropy parameter. Note that for $r=1$  the energy reduces to that of a neo-Hookean material. 
As most nematic rubbers are nearly incompressible \cite{WT03},
we prescribe $W^{e}$ to be a finite  only for volume-preserving deformations.   
We can substitute for $\ell$ and write
\begin{eqnarray}\label{eq:We}
W^{e}(F,n)= \begin{cases} 
\displaystyle{\frac{\mu}{2}\Bigl(r^{1/3} \Bigl(|F|^2-\frac{r-1}{r}|F^Tn|^2\Bigr)-3\Bigr)} &  \textrm{if } \det F=1, |n| = 1, \\
+\infty & \textrm{else} .
\end{cases}
\end{eqnarray}

For future use, we define a purely elastic energy by taking the infimum over directors.  Following DeSimone and Dolzmann \cite{DD00}, 
\begin{eqnarray}\label{eq:W3D}
W_{3D}(F):=\inf_{n \in \mathbb{S}^2}W^{e}(F,n)= \begin{cases}
\displaystyle{W_0(F)}&  \textrm{if }\det F=1, \\
+\infty & \textrm{else},
\end{cases} 
\end{eqnarray}
where
\begin{eqnarray}\label{eq:W0}
W_0(F)=\frac{\mu}{2}\Bigl(r^{1/3} \Bigl(|F|^2-\frac{r-1}{r}\lambda_{M}^2(F)\Bigr)-3\Bigr).
\end{eqnarray}
Here $\lambda_{M}(F)$ is the largest eigenvalue of $(F^TF)^{1/2}$.
Energy density $W_{3D}(F)$ is not quasiconvex and the quasiconvex envelope $W^{qc}_{3D}$ has been computed in \cite{DD00}.

Following Oseen, Zocker and Frank, (see for example, \cite{DP93}), we take the elasticity of the nematic mesogens to be of the form
\begin{eqnarray}\label{eq:Wn}
W^{n} = \frac{1}{2} \kappa_1 (\divg n )^2 + \frac{1}{2} \kappa_2 ( n \cdot \curl n )^2 + \frac{1}{2} \kappa_3 (n \times \curl n )^2
\end{eqnarray}
where $\divg n$ and $\curl n$ are the spatial divergence and curl of the director respectively, and $\kappa_1>0, \kappa_2>0, \kappa_3>0$ are known as the splay, twist and bend moduli respectively.   Notice that this is a non-negative quadratic form in $n$ and $\nabla_y n$.  It turns out that these moduli are very close to each other and one can introduce an equal modulus approximation 
\begin{eqnarray}
W^{n}_{eq} = \frac{\kappa}{2}  |\nabla_y n|^2 =  \frac{\kappa}{2} |(\nabla n) F^{-1}|^2
\end{eqnarray}
\noindent where the second equality holds formally. 
Importantly, from a mathematical point of view, any given $W^n$ of the form (\ref{eq:Wn}) can be bounded from above and below by equal moduli approximations, and therefore all the results we prove for the equal modulus approximation hold for the more general form.    Finally, since we assume incompressibility or $\det F = 1$, $F^{-1} = \adj F$ so that 
\begin{eqnarray}\label{eq:Wneq}
W^{n}_{eq} = \frac{\kappa}{2} |(\nabla n) (\adj F) |^2.
\end{eqnarray}

Putting these together, the Helmholtz free energy of our nematic elastomer, $\mathcal{F}:H^1(\Omega,\mathbb{R}^3) \times W^{1,1}(\Omega,\mathbb{S}^2) \rightarrow \mathbb{R} \cup \{ +\infty\}$, is given by  
\begin{equation}\label{eq:HFE}
\mathcal{F}(y,n) :=\begin{dcases}
 \int_\Omega \left( W^{e} (\nabla y, n) + \frac{\kappa}{2} | (\nabla n) (\adj \nabla y)|^2 \right) dx &\text{ if } (y,n) \in \mathcal{A}, \\
 +\infty &\text{ else } ,
 \end{dcases}
\end{equation}
\noindent where for definiteness, the set of admissible fields is
\begin{align*}
\mathcal{A} := \{ (y,n) \in H^1(\Omega,\mathbb{R}^3) \times W^{1,1}(\Omega,\mathbb{S}^2): (\nabla n)(\adj \nabla y) \in L^2(\Omega,\mathbb{R}^{3\times3}) \}.  
\end{align*}
\noindent The ambient space for deformations, $y \in H^1(\Omega,\mathbb{R}^3)$, is optimal since for $F \in \mathbb{R}^{3\times3}$ satisfying $\det F = 1$, $W^e$ satisfies the growth and coercivity
\begin{align}\label{eq:WeG}
\frac{1}{c}|F|^2 - c  \leq W^e(F,n) \leq c(|F|^2 + 1) 
\end{align}
\noindent independent of $n \in \mathbb{S}^2$.  Here $c \geq 1$ depends on $r$ and $\mu$.  The ambient space for the director field, $n \in W^{1,1}(\Omega,\mathbb{S}^2)$, may not be optimal.  Nevertheless, consider the following: 

\begin{remark}
Fonseca and Gangbo \cite{FG95} showed the lower-semicontinuity of 
\begin{align*}
\int_{\Omega}|(\nabla n)(\nabla y)^{-1}|^2 dx 
\end{align*}
in the space 
\begin{align*}
\beta_{p,q} = \{ (n,y) \in W^{1,p}(\Omega,\mathbb{S}^2) \times W^{1,q}(\Omega,\mathbb{R}^3) : \det \nabla y (x) = 1 \text{ a.e. } x \in \Omega\}, 
\end{align*}
where $2 <  p < + \infty$ and $4 <  q \leq +\infty$ such that $(1/p)+(2/q) = (1/2)$, (see Theorem 4.1, \cite{FG95}).    The existence of minimizers in $\beta_{p,q}$ follows from this result.  However, notice that this is more regularity than we assume.  In fact, the existence of minimizers in ${\mathcal A}$ is not clear.   However, this does not affect $\Gamma-$convergence or the membrane limit.  

Barchiesi and DeSimone \cite{BD13} showed  well-posedness for an energy similar to (\ref{eq:HFE}) where $n$ was taken as a mapping from the deformed configuration $y(\Omega)$, and it was assumed $W^e(F,n) \geq c(|F|^3 - 1)$. 
\end{remark}

\section{Membrane Theory}\label{MemTheory}

In this section, we derive a theory for nematic elastomer membranes whose three dimensional free energy satisfies (\ref{eq:HFE}).

\subsection{Framework}

We consider a nematic elastomer membrane of small thickness $h$ which has a flat stress-free isotropic reference configuration $\Omega_h := \{ (x',x_3) \in \mathbb{R}^3 : x' \in \omega, x_3 \in (-h/2,h/2)\}$.  We assume $\omega$ is a bounded Lipschitz domain in $\mathbb{R}^2$.  Let $\tilde{y} :\Omega_h \rightarrow \mathbb{R}^3$ describe the deformation and $\tilde{n} : \Omega_h \rightarrow \mathbb{S}^2$ describe the director field so that $\mathcal{F}^h(\tilde{y},\tilde{n})$ is the Helmholtz free energy in (\ref{eq:HFE}) now parameterized by the thickness of the membrane in its reference configuration.  We assume $\kappa/2 = \kappa_h$, $\kappa_h \geq 0$ and $\kappa_h \rightarrow 0$ as $h \rightarrow 0$.  

To take the limit as $h \rightarrow 0$, we follow the theory of $\Gamma$-convergence in a topological space endowed with the weak topology.  The general theory can be found in \cite{B02} and \cite{D93}.  In order to deal with sequences on a fixed domain, we change variables via
\begin{align*}
z' = (z_1,z_2) = (x_1,x_2) = x', \;\;\;\; z_3 = \frac{1}{h} x_3, \;\;\;\; x \in \Omega_h
\end{align*} 
\noindent and set $\Omega := \omega \times (-1/2,1/2)$.  To each deformation $\tilde{y} : \Omega_h \rightarrow \mathbb{R}^3$ and director field $\tilde{n}:\Omega_h \rightarrow \mathbb{S}^2$ , we associate respectively a deformation $y: \Omega \rightarrow \mathbb{R}^3$ and director field $n: \Omega \rightarrow \mathbb{S}^2$ such that 
\begin{align}\label{eq:ChVar2}
y(z(x)) = \tilde{y}(x) \;\;\;\; \text{ and } \;\;\;\; n(z(x)) = \tilde{n}(x), \;\;\;\; x \in \Omega_h.
\end{align}
\noindent We set $\mathcal{\tilde{I}}^h(y,n) := \mathcal{F}^h(\tilde{y},\tilde{n}) /h $, and following the change of  variables above observe
\begin{align}\label{eq:UBIh}
\mathcal{\tilde{I}}^h(y,n) = \begin{dcases}
\int_{\Omega} \left(W^e(\nabla_h y,n) + \frac{\kappa_h}{h^2} |(\nabla n )(\adj \nabla y)|^2 \right) dz &\text{ if } (y,n) \in \mathcal{A}, \\
+\infty &\text{ else }
\end{dcases}
\end{align}  
where $\nabla_h y = (\nabla' y|(1/h)\partial_3 y)$ with $\nabla'$ the in-plane gradient.  We also use the identity $(\nabla_h n)(\adj \nabla_h y) = (1/h) (\nabla n) (\adj \nabla y)$. 

Finally, we take our membrane theory to be the $\Gamma$-limit as $h \rightarrow 0$ of the functional defined on $H^1(\Omega,\mathbb{R}^3)$,
\begin{align}\label{eq:Ih}
\mathcal{I}^h(y) := \inf_{n \in W^{1,1}(\Omega,\mathbb{S}^2)} \mathcal{\tilde{I}}^h(y,n).  
\end{align}


\subsection{The Membrane Limit}

\begin{theorem}\label{1407251605}
Let $\mathcal{I}^h$ be as in (\ref{eq:Ih}) with $\kappa_h \geq 0$ and $\kappa_h \rightarrow 0$ as $h \rightarrow 0$.  Then in the weak topology of $H^1(\Omega,\mathbb{R}^3)$, $\mathcal{I}^h$ is equicoercive and $\Gamma$-converges to 
\begin{align}\label{eq:J}
\mathcal{J}(y) = \begin{dcases}
\int_{\omega} W_{2D}^{qc}( \nabla ' y) dz' &\text{ if } \partial_3 y = 0 \text{ a.e.}, \\
+\infty &\text{ otherwise }.
\end{dcases}
\end{align}
\noindent Here
\begin{align}\label{eq:W2DTh}
W_{2D}(\tilde{F}) := \inf_{c \in \mathbb{R}^3} W_{3D}(\tilde{F}|c) 
\end{align}
\noindent for $W_{3D}$ given in (\ref{eq:W3D}) and $W_{2D}^{qc}$ is the quasiconvex envelope of $W_{2D}$,
\begin{align*}
W_{2D}^{qc}(\tilde{F}) = \inf \left\{ \int_{(0,1)^2} W_{2D}(\tilde{F} + \nabla ' \phi) dz': \phi \in W^{1,\infty}_0((0,1)^2,\mathbb{R}^3) \right\}.
\end{align*}
\noindent Equivalently:
\begin{enumerate}[(i)]
\item for every sequence $\{y_h\} \subset H^1(\Omega,\mathbb{R}^3)$ such that $\mathcal{I}^h(y_h) \leq C < +\infty$, there exists a $y \in H^1(\Omega,\mathbb{R}^3)$ independent of $z_3$ such that up to a subsequence 
\begin{align*}
y_h - \fint_{\Omega} y_h dz \rightharpoonup y \;\;\;\; \text{ in } H^1(\Omega,\mathbb{R}^3);
\end{align*}
\item for every $\{ y_h \} \subset H^1(\Omega,\mathbb{R}^3)$ such that $y_h \rightharpoonup y$ in $H^1(\Omega,\mathbb{R}^3)$, 
\begin{align*}
\liminf_{h \rightarrow 0} \mathcal{I}^h(y_h) \geq \mathcal{J}(y);
\end{align*}
\item for any $y \in H^1(\Omega,\mathbb{R}^3)$, there exists a sequence $\{ y_h \} \subset H^1(\Omega,\mathbb{R}^3)$ such that $y_h \rightharpoonup y$ in $H^1(\Omega,\mathbb{R}^3)$ and 
\begin{align*}
\limsup_{h \rightarrow 0} \mathcal{I}^h(y_h) \leq \mathcal{J}(y).  
\end{align*}
\end{enumerate} 
\end{theorem}

The result for the case $\kappa_h = 0$ was provded by Conti and Dolzmann \cite{CD06} (Theorem 3.1 there).

\begin{theorem}[Conti and Dolzmann \cite{CD06}]\label{ContiDolzmannResult}
In the weak topology of $H^1(\Omega,\mathbb{R}^3)$, the functional
\begin{align*}
\mathcal{I}_e^h(y) := \int_{\Omega} W_{3D}(\nabla_h y) dz 
\end{align*}
\noindent is equicoercive and $\Gamma$-converges to $\mathcal{J}$ given in (\ref{eq:J}).  
\end{theorem}

\begin{remark}
A different dimension reduction theory for hyperelastic incompressible materials was developed by Trabelsi \cite{Trabelsi05},\cite{Trabelsi06} under similar assumptions. Trabelsi shows that the membrane energy density (integrand of ${\mathcal J}$)
is given by $((W_{2D})^{rc})^{qc}$. 
From the proof of Theorem \ref{1407231520} below, it follows that $W_{2D}^{rc}=W_{2D}^{qc}$ (and hence $(W_{2D}^{rc})^{qc}=W_{2D}^{qc}$).  Thus the two limits agree.
\end{remark}

\begin{remark}\label{1407231441}
We remark on some general properties of the purely elastic portion of our nematic elastomer energy density. $W_0:\mathbb{R}^{3\times3} \rightarrow \mathbb{R}$ in (\ref{eq:W0}) is Lipschitz continuous, $W_0 + 3\mu/2$ is non-negative, and there exists a constant $c$ such that
\begin{align}\label{eq:W0C}
\frac{1}{c}|F|^2 - c \leq W_0(F) \leq c(|F|^2 + 1).
\end{align}
The energy $W_{2D}$ in (\ref{eq:W2DTh}) is given by 
\begin{align}\label{eq:W2D1}
W_{2D}(\tilde{F}) = \begin{dcases}
\min_{c \in \mathbb{R}^3} W_{3D}(\tilde{F}|c) &\text{ if } \rank \tilde{F} = 2,\\
+\infty &\text{ else },
\end{dcases}
\end{align}
\noindent and satisfies 
\begin{align}\label{eq:W2DG}
\frac{1}{c} \left( |\tilde{F}|^2 + \frac{1}{ \delta(\tilde{F})^2} \right) - c \leq W_{2D}(\tilde{F}) \leq c\left( |\tilde{F}|^2 + \frac{1}{\delta(\tilde{F})^2} + 1 \right)
\end{align}
\noindent with $\delta(\tilde{F}) = |\adj(\tilde{F})|$.  The \textit{effective} energy density $W_{2D}^{qc}$ is quasiconvex, Lipschitz continuous on bounded sets and its definition does not depend on the choice of the domain $\omega$, as long as it is open, bounded and $|\partial\omega|=0$.
Furthermore, there exists (Lemma 3.1, \cite{CD06}) a constant $c'$ such that 
\begin{eqnarray}\label{eq:W2DQCBound}
\frac{1}{c'}|\tilde F|^2-c'\leq W_{2D}^{qc}(\tilde{F})\leq c'|\tilde{F}|^2+c'.
\end{eqnarray}
\end{remark} 
\bigskip

\noindent {\it Proof of Theorem \ref{1407251605}.}
Note that trivially, 
\begin{align}\label{eq:IehBound}
\mathcal{I}^h(y) \geq  \int_{\Omega} \inf_{n \in \mathbb{S}^2} W^e(\nabla_h y,n) dz = \mathcal{I}_e^h(y).  
\end{align}
Therefore, the compactness  and lower bound (Properties (i) and (ii) in Theorem \ref{1407251605}) follow
from Theorem \ref{ContiDolzmannResult}.  It remains to show Property (iii).  This is done in Proposition \ref{1407251610}.
\endproof

Before we proceed, we note that the fact that the $\Gamma-$limit is independent of $\kappa_h/h$ is similar to the following result of Shu \cite{Shu00}.  He also provides some heuristic insight.  Since the membrane limit optimizes the energy density over the third column of the deformation gradient, there is little to be gained by oscillations  parallel to the thickness.  Consequently, penalizing these oscillations with $\kappa_h$ does not affect the $\Gamma-$limit.

\begin{theorem}[Shu \cite{Shu00}] \label{YCShuResult}
Let $\kappa_h \rightarrow 0$ as $h \rightarrow 0$, and $W: \Omega \rightarrow {\mathbb R}$ be continuous and bounded from above and below by $|F|^p \pm c$ respectively for some $c$.  Then, in the weak topology of $W^{1,p}$, the functional
\begin{align*}
\int_{\Omega}\{ \kappa_h | \nabla_h \nabla_h y|^2 + W( \nabla_h y) \}dz,
\end{align*}
$\Gamma-$converges to 
\begin{align*}
\begin{cases} \int_{\omega}\{ W_{2D}^{qc}( \nabla' y) \}dz'  & if \  \partial_3 y = 0 \ a.e., \\
\infty & else ,
\end{cases}
\end{align*}
where $W_{2D}(\tilde{F}) := \inf_{c \in \mathbb{R}^3} W_{3D}(\tilde{F}|c)$.
\end{theorem}


\subsection{Construction of Recovery Sequence}\label{RecSeq}

It remains to construct a recovery sequence to prove $\mathcal{J}$ is the $\Gamma$-limit to $\mathcal{I}^h$.    

\begin{proposition}\label{1407251610}
For every $y \in H^1(\Omega,\mathbb{R}^3)$ independent of $z_3$, there exists a sequence $\{ (y_h,n_h) \} \subset 
C^{\infty}(  \bar{\Omega},\mathbb{R}^3) \times 
C^{1}( \bar{\Omega},\mathbb{S}^2)$ such that $y_h \rightharpoonup y $ in $H^1(\Omega,\mathbb{R}^3)$ and 
\begin{align}\label{eq:PropEqn}
\limsup_{h \rightarrow 0} \mathcal{\tilde{I}}^h(y_h,n_h)  \leq \mathcal{J}(y). 
\end{align}
\end{proposition}

Our construction also draws heavily from Conti and Dolzmann \cite{CD06}.  The main difference is that we need additional regularity for our recovery sequence $n_h$.   We summarize the Conti-Dolzmann construction in two lemmas.  The first lemma regards the construction of a sequence  to go from the energy density  $W_{2D}$ to $W_{2D}^{qc}$ on $\omega$.  For our analysis, the important observation is that in the limit the deformation gradient is constant on an increasingly large subset of $\omega$. The second lemma regards the extension of smooth maps on $\omega$ to incompressible deformations on $\Omega_h$.  

\begin{lemma}[S. Conti and G. Dolzmann \cite{CD06}]\label{1407251613}
For any $y \in H^1(\omega,\mathbb{R}^3)$, there exists a sequence $\{y_j \} \subset C^{\infty}(\bar{\omega},\mathbb{R}^3)$ such that $\rank \nabla y_j = 2$ everywhere, $y_j \rightharpoonup y$ in $H^1(\omega,\mathbb{R}^3)$ as $j \rightarrow \infty$, and 
\begin{align}\label{eq:seqWqc}
\limsup_{j \rightarrow \infty} \int_{\omega} W_{2D}(\nabla' y_j) dx' \leq \int_{\omega} W_{2D}^{qc}(\nabla' y) dx'.
\end{align}
\noindent Moreover, the sequence has the following properties:
\begin{enumerate}[(i)]
\item for each $j \in \mathbb{N}$, $y_j$ is defined on a triangulation $\mathcal{T}^{\;j}$ of $\omega$ which is the set of at most countably many disjoint open triangle $T_i^{\;j}$ whose union up to a null set is equal to $\omega$, and $\Gamma_j$ is the jump set given by 
\begin{align*}
\Gamma_j := \partial \omega \cup \bigcup_{i} \partial T_i^{\;j};
\end{align*}
\item there is a sequence of boundary layers $\{ \eta_j \}$ such that $\eta_j > 0$ and $\eta_j \rightarrow 0 $ as $j \rightarrow \infty$, and the set $\Gamma_{\eta_j}$ is defined to be 
\begin{align*}
\Gamma_{\eta_j} := \{ x' \in \omega: \dist(x', \Gamma_j ) <  \eta_j \};
\end{align*}
\item if $T_i^{\;j} \setminus \Gamma_{\eta_j} $ is nonempty, then $\nabla' y_j$ is a constant on this set and we set 
\begin{align}\label{eq:constF}
\tilde{F}_{i}^{\;j} := \nabla y_j (x'), \;\;\;\; x' \in T_{i}^{\;j} \setminus \Gamma_{ \eta_j};
\end{align}
\item   $\adj \nabla' y_j$ is bounded away from zero in the sense that for some $\epsilon_j > 0$ sufficiently small, the inequality 
\begin{align*}
|\adj \nabla' y_j| \geq \epsilon_j > 0,
\end{align*} 
holds everywhere.  
\end{enumerate}
\end{lemma}
  
\begin{lemma}[Conti and Dolzmann \cite{CD06}]\label{1407251618}
Let $w, \nu \in C^{\infty}(\bar{\omega}, \mathbb{R}^3)$ satisfy 
\begin{align*}
\det(\nabla' w|\nu) = 1 \;\;\;\; \text{ in } \omega.   
\end{align*}
\noindent Then there exists an $h_0 >0$ and an extension $v \in C^{\infty}(\bar{\omega} \times (-h_0,h_0), \mathbb{R}^3)$ such that $v(x',0) = w(x')$ and $\det \nabla v = 1$ everywhere.  Moreover, for all $x_3 \in (-h_0,h_0)$ the pointwise bound 
\begin{align*}
|\nabla v(x) - (\nabla ' w| \nu)(x') | \leq C |x_3|
\end{align*}
\noindent holds, where $C$ can depend on $w$ and $\nu$.  
\end{lemma}

We construct a recovery sequence and thereby prove Proposition \ref{1407251610} in 4 parts.  In Part 1, we take a sequence of smooth maps $y_j$ as in Lemma \ref{1407251613} and show that we can construct a sequence of smooth vector fields $c_j$  such that $\det (\nabla' y_j | c_j) = 1$ in $\omega$.  In Part 2, we use Lemma \ref{1407251618} to extend $y_j$ appropriately to a deformation on $\Omega$, i.e $y_j^{\;h}$.  In Part 3, we construct a sequence of $C^1$ director fields $n_j^{\;h}$ on $\Omega$ which enables passage from $W^e$ to $W_{2D}$.  Finally, in Part 4 we show that we can take an appropriate diagonal sequence $h_j \rightarrow 0$ as $j \rightarrow \infty$ which proves Proposition \ref{1407251610}. \\

\noindent{\it Proof of Proposition \ref{1407251610}.}  Let $y \in H^1(\Omega,\mathbb{R}^3)$ independent of $z_3$.  Then $y$ is bounded in $H^1(\omega,\mathbb{R}^3)$ (with abuse of notation). By Lemma \ref{1407251613}, we find a sequence $\{ y_j \} \subset C^{\infty}(\bar{\omega},\mathbb{R}^3)$ such that $\rank \nabla' y_j = 2$ everywhere, $y_j \rightharpoonup y $ in $H^1(\omega,\mathbb{R}^3)$, the energy is bounded in the sense of (\ref{eq:seqWqc}), and the sequence satisfies properties (i)-(iv) from the lemma.

\noindent Part 1.   We define the smooth vector field $c_j$ on the triangulation $\mathcal{T}^{\;j}$ for $y_j$ in Lemma \ref{1407251613} (i).  On each nonempty $T_i^{\;j}\setminus \Gamma_{\eta_j}$ there exists a constant $\tilde{F}_i^{\;j}$ defined in Lemma \ref{1407251613} (iii), and it is full rank.  Then by (\ref{eq:W2D1}), $W_{3D}(\tilde{F}_i^{\;j}|c)$ has a minimizer for $c \in \mathbb{R}^3$.  Motivated by this observation, we let
\begin{align}\label{eq:cj1}
c_i^{\;j} := \arg \min_{c \in \mathbb{R}^3}  W_{3D}(\tilde{F}_i^{\;j}|c),
\end{align}
\noindent which via (\ref{eq:W3D}) implies
\begin{align}\label{eq:det1}
\det (\tilde{F}_i^{\;j}|c_i^{\;j}) = 1.
\end{align}

Consider the vector field,
\begin{align}\label{eq:cj2}
c_0^{\;j} := \frac{\adj \nabla' y_j}{|\adj \nabla' y_j|^2} .
\end{align}
\noindent  This is well-defined given Lemma \ref{1407251613} (iv).  Moreover, since $y_j$ is smooth, $c_0^{\;j} \in C^{\infty}(\bar{\omega},\mathbb{R}^3)$.  Further, since $\det(\tilde{F} | c) = (\adj \tilde{F})^T c $, we have 
\begin{align}\label{eq:det2}
\det(\nabla ' y_j |c_0^{\;j}) = 1 \;\;\;\; \text{ in } \omega.
\end{align}

Let $c_j \in C^{\infty}(\bar{\omega}, \mathbb{R}^3)$ be given by 
\begin{align}\label{eq:cj}
c_j := \begin{dcases}
c_0^{\;j} + \psi_i \left( c_i^{\;j} - c_0^{\;j} \right) &\text{ on each } T_i^{\;j}\setminus \Gamma_{\eta_j} \text{ with nonempty open subsets}, \\
c_0^{\;j} &\text{ otherwise on } \omega.
\end{dcases}
\end{align}
\noindent Here $\psi_i \in C^{\infty}_0(T_i^{\;j}\setminus \Gamma_{\eta_j}, [0,1])$ is a cutoff function which equals 1 at least on the entirety of the subset $T_i^{\;j} \setminus \Gamma_{2\eta_j}$.  Notice when $c_j = c_0^{\;j}$, the determinant constraint is satisfied trivially by (\ref{eq:det2}).  Conversely, combining (\ref{eq:det1}) and (\ref{eq:det2}),
\begin{align*}
\det (\nabla' y_j | c_j) &= (\adj \nabla' y_j)^T  \left(c_0^{\;j} + \psi_i \left( c_i^{\;j} - c_0^{\;j} \right)\right) \nonumber\\ &= \det(\nabla ' y_j | c_0^{\;j}) + \psi_i \left( \det(\tilde{F}_{i}^{\;j} |c_i^{\;j}) - \det(\nabla' y_j | c_0^{\;j} ) \right) \nonumber \\
&= 1 \;\;\;\; \text{ on each } T_i^{\;j} \setminus \Gamma_{\eta_j} \text{ with nonempty open subsets },  
\end{align*}
\noindent since $\nabla' y_j = \tilde{F}_{i}^{\;j}$ on this set.  We then conclude $\det(\nabla' y_j |c_j) = 1$ in $\omega$, and this completes Part 1.  

\noindent Part 2.  Fix $j \in \mathbb{N}$.  From Part 1 we have $y_j,c_j \in C^{\infty}(\bar{\omega}, \mathbb{R}^3)$ satisfying $\det( \nabla' y_j |c_j) = 1$ in $\omega$.  Hence, there exists an $h_0^{\;j} > 0$ and a $v \in C^{\infty}(\bar{\omega}\times (-h_0^{\;j},h_0^{\;j}))$ such that the properties of Lemma \ref{1407251618} hold replacing $w$ with $y_j$ and $\nu$ with $c_j$.  Let $h \in (0,h_0^{\;j})$ and $\tilde{y}_j^{\;h} \in 
C^{\infty}(  \bar{\Omega}_h,\mathbb{R}^3)$ be the restriction of $v$ to $\Omega_h$.  Further, let $y_j^{\;h} \in 
C^{\infty}( \bar{\Omega},\mathbb{R}^3)$ be associated to $\tilde{y}_j^{\;h}$ using (\ref{eq:ChVar2}).  From Lemma \ref{1407251618}, we conclude $y_j^{\;h}(z',0) = y_j(z')$, $\det \nabla_h y_j^{\;h}(z) = 1$ and 
\begin{align}\label{eq:bGrad}
|\nabla _h y_j^{\;h}(z) - (\nabla ' y_j | c_j)(z')| \leq C_j h |z_3| \leq C_j h, \;\;\;\; z \in \Omega.
\end{align}
\noindent Here $C_j$ is a constant depending on $y_j$ and $c_j$, and the second inequality above follows since $z_3 \in (-1/2,1/2)$.   From these properties we conclude as $h \rightarrow 0$,
\begin{align}\label{eq:yjhCon}
y_j^{\;h} \rightarrow y_j \text{ in } H^1(\Omega,\mathbb{R}^3)\;\;\;\; \text{ and } \;\;\;\; \frac{1}{h} \partial_3 y_j^{\;h} \rightarrow c_j \text{ in } L^{2}(\Omega,\mathbb{R}^3).
\end{align}
\noindent This concludes Part 2.  

\noindent Part 3.  As in Part 2, we keep $j \in \mathbb{N}$ fixed.  From Lemma \ref{1407251613} (i) we have that $\bigcup_i T_i^{\;j} = \omega$ (up to a set of zero measure), though this union can be countably infinite.  From herein, we choose a finite collection of $N(j)$ triangles so that
\begin{align}\label{eq:Nj}
\int_{\omega\setminus \cup_{i =1}^{N(j)}T_i^{\;j}} \left\{ W_{2D}(\nabla' y_j) + 1 \right\} dz' \leq \frac{1}{j}.  
\end{align}
\noindent Then for each of the $N(j)$ triangles for which the set $T_{i}^{\;j} \setminus \Gamma_{\eta_j}$ is nonempty, let 
\begin{align}\label{eq:nij}
n_i^{\;j} := \arg \min_{n \in \mathbb{S}^2} W^e(\tilde{F}_i^{\;j}|c_i^{\;j},n) .
\end{align}
\noindent Further, let $q_j$ be the piecewise constant function on $\mathbb{R}^2$ given by 
\begin{align}\label{eq:prenj}
q_j(z') := \begin{dcases} 
n_i^{\;j} &\text{ if } i \in \{1,\ldots, N(j)\}, \;T_{i}^{\;j} \setminus \Gamma_{\eta_j} \text{ is nonempty}, \text{ and } z' \in T_i^{\;j},  \\
q &\text{ otherwise in } \mathbb{R}^2.  
\end{dcases} 
\end{align}
\noindent Here $q$ is a fixed vector in $\mathbb{S}^2$.  Then $q_j$ maps to $\mathbb{S}^2$, but it is not in $C^1$.  To correct this, we employ the approach used by DeSimone in \cite{DS93} (see Assertion 1).  

Observe by construction the range of $q_j$ is finite.  Hence, there exists an $s_j \in \mathbb{S}^2$ and a closed ball $B_{\epsilon}(s_j)$ of radius $\epsilon > 0$ centered at $s_j$ such that $(\range q_j) \cap B_{\epsilon}(s_j) = \emptyset$.  Then the stereographic projection $\pi_{s_j}$  with the projection point as $s_j$ maps the range of $q_j$ to a bounded subset of $\mathbb{R}^2$.  Let $\psi_{\eta_j}$ be a standard mollifier with $\eta_j$ as in Lemma \ref{1407251613} (ii), and consider the composition 
\begin{align*}
\tilde{n}_j := \pi_{s_j}^{-1} \circ \left( \psi_{\eta_j} \ast \left( \pi_{s_j} \circ q_j\right) \right).
\end{align*}  
\noindent  This composition is well-defined since the range of $q_j$ is outside a neighborhood of the projection point $s_j$.  Further, $\tilde{n}_j$ maps to $\mathbb{S}^2$ using the definition of the inverse of the stereographic projection.  Moreover, $\pi_{s_j}^{-1}$ is differentiable and its argument $\psi_{\eta_j} \ast ( \pi_{s_j} \circ q_j )$ is smooth.  Hence, $\tilde{n}_j \in C^1(\mathbb{R}^2,\mathbb{S}^2)$.  

Let $n_j \in 
C^1( \bar{\omega},\mathbb{S}^2)$ be the restriction of $\tilde{n}_j$ to the closure of $\omega$.  Further, let $n_j^{\;h} \in 
C^1( \bar{\Omega},\mathbb{S}^2)$ be the extension of $n_j$ to $\Omega$ via $n_j^{\;h} (z) := n_j(z')$ for each $z \in \Omega$.   As a final remark for this part, observe for $i \in \{1, \ldots, N(j)\} $ and $z' \in T_i^{\;j} \setminus \Gamma_{2\eta_j}$, 
\begin{align}\label{eq:constn}
n_j^{\;h}(z) = n_j(z') &=  \pi_{s_j}^{-1} \circ \left( \psi_{\eta_j} \ast \left( \pi_{s_j} \circ q_j\right) \right) (z') \nonumber \\
&= \pi_{s_j}^{-1} \circ \left( \int_{\mathbb{R}^2} \psi_{\eta_j}(z'- \xi) (\pi_{s_j} \circ q_j)(\xi) d\xi \right) \nonumber \\
&= \pi_{s_j}^{-1} \circ \left( (\pi_{s_j} \circ q_j) \int_{B_{\eta_j}(z')} \psi_{\eta_j}(z'-\xi) d\xi\right) \nonumber \\
&= \pi_{s_j}^{-1} \circ (\pi_{s_j} \circ q_j) = n_i^{\;j},  
\end{align}
\noindent  since $B_{\eta_j}(z') \cap \partial T_i^{\;j} = \emptyset$ and so $q_j$ is constant on $B_{\eta_j}(z')$, see (\ref{eq:prenj}).  This completes Part 3.  

\noindent Part 4.  From Parts 1-3, we have for each $j \in \mathbb{N}$ the functions $y_j^{\;h} \in C^{\infty}( \bar{\Omega},\mathbb{R}^3)$ and $n_j^{\;h} \in 
C^1( \bar{\Omega}, \mathbb{S}^2)$ parameterized by $h \in (0,h_0^{\;j})$.  It remains to bound the functional $\mathcal{I}^h$ appropriately and take the $\limsup$.  For the bounding arguments, $C$ shall refer to positive constant independent of $h$ and $j$ which may change from line to line. From (\ref{eq:We1}), when $W^e$ is finite, it satisfies a Lipschitz condition 
\begin{align*}
|W^e(F,n)- W^e(G,n)| &\leq |\ell^{-1/2}|^2\left( |F| + |G| \right)|F - G| \nonumber \\
&\leq C \left( |F| + |G| \right) |F - G| .
\end{align*}
\noindent As asserted above, $|\ell^{-1/2}|$ is uniformly bounded for $n \in \mathbb{S}^2$.  Then since for every $z \in \Omega$, $\det (\nabla_h y_j^{\;h})(z) = 1$, $\det (\nabla' y_j |c_j )(z') = 1$ and $n_j^{\;h}(z) = n_j(z') \in \mathbb{S}^2$, 
\begin{align}\label{eq:UB1}
\int_{\Omega} W^e(\nabla_h y_j^{\;h}, n_j^{\;h}) dz \nonumber 
&\leq \int_{\omega} W^e(\nabla ' y_j|c_j,n_j) dz' + \int_\Omega| W^e(\nabla_h y_j^{\;h},n_j^{\;h}) - W^e(\nabla' y_j| c_j, n_j^{\;h}) | dz \nonumber \\
&\leq \int_{\omega} W^e(\nabla ' y_j|c_j,n_j) dz' + E_{h,j}^1.
\end{align} 
\noindent Here $E_{h,j}^1$ is the estimate obtained from the Lipschitz condition and an application of H\"older's inequality,
\begin{align}\label{eq:E1hj}
E_{h,j}^1 := C\left( \|\nabla_h y_j^{\;h} \|_{L^2(\Omega,\mathbb{R}^3)} +\| (\nabla ' y_j| c_j) \|_{L^2(\Omega,\mathbb{R}^3)} \right) \| \nabla_h y_j^{\;h} - (\nabla' y_j| c_j) \|_{L^2(\Omega,\mathbb{R}^3)}.
\end{align}

We now focus on the first term in the upper bound (\ref{eq:UB1}).  For $i \in \{1,...,N(j)\}$ and  $z' \in T_i^{\;j} \setminus \Gamma_{2\eta_j}$, observe 
\begin{align*}
W^e(\nabla' y_j(z')|c_j(z'), n_j(z'))  &= W^e(\tilde{F}_i^{\;j} | c_i^{\;j}, n_i^{\;j}) \;\;\;\;\;\;\;\;\;\;\;\;\;\;\;\;\;\;\;\;\;\;\; \text{ by } (\ref{eq:constF}), (\ref{eq:cj}) \text{ and } (\ref{eq:constn}); \nonumber \\
&= \min_{n \in \mathbb{S}^2} W^e(\tilde{F}_i^{\;j}|c_i^{\;j}, n) \;\;\;\;\;\;\;\;\;\;\;\;\;\;\;\;\;\;\text{ by } (\ref{eq:nij}); \nonumber  \\
&= \min_{c \in \mathbb{R}^3} W_{3D}( \tilde{F}_i^{\;j} | c)\;\;\;\;\;\;\;\;\;\;\;\;\;\;\;\;\;\;\;\;\;\;\text{ by } (\ref{eq:W3D}) \text{ and }(\ref{eq:cj1}); \nonumber \\
&= W_{2D}(\tilde{F}_i^{\;j} ) = W_{2D}(\nabla' y_j(z')) \;\;\;\; \text{ by } (\ref{eq:W2D1}).
\end{align*} 
\noindent Then,
\begin{align}\label{eq:UB2}
\int_{\omega} W^e(\nabla ' y_j|c_j,n_j) dz' &\leq \int_{(\cup_{i = 1}^{N(j)}T_{i}^{\;j})\setminus \Gamma_{2\eta_j}} W_{2D}(\nabla ' y_j) dz' + \int_{\omega \setminus \cup_{i = 1}^{N(j)}T_{i}^{\;j}} W^e(\nabla ' y_j | c_j,n_j) dz' \nonumber \\
&\;\;\;\; + \int_{\Gamma_{2\eta_j}} W^e(\nabla ' y_j | c_j,n_j) dz',
\end{align}
\noindent using our result for $W^e$ and since each integrand is nonnegative. 

We bound $W^e(\nabla' y_j| c_j,n_j)$ in (\ref{eq:UB2}).  To obtain this bound notice $\; |c_0^{\;j}|^2 = 1/|\adj \nabla' y_j|^2$ from (\ref{eq:cj2}).  Further, using the coercivity condition of $W_0$ in (\ref{eq:W0C}), the definition of $c_i^{\;j}$ in (\ref{eq:cj1}), and the growth in (\ref{eq:W2DG}),
\begin{align*}
|c_i^{\;j}|^2 \leq W_{0}(\tilde{F}_i^{\;j}| c_i^{\;j}) = W_{2D}(\tilde{F}_{i}^{\;j}) \leq c\left(|\tilde{F}_i^{\;j}|^2 + \frac{1}{|\adj \tilde{F}_i^{\;j}|^2} + 1\right). 
\end{align*}
\noindent 
Following these observations, we notice on the sets $T_i^j \setminus \Gamma_{\eta_j}$,  $\nabla y_j = \tilde{F}_i^{\;j}$ by definition (see Lemma \ref{1407251613} (iii)) and therefore,
\begin{align*}
|c_j|^2 \leq 2(|c_0^{\;j}|^2 + |c_i^{\;j}|^2) \leq C\left( |\nabla' y_j|^2 + \frac{1}{|\adj \nabla' y_j|^2} + 1\right)
\end{align*}
since $c_j$ is as in (\ref{eq:cj}).  On the exceptional sets, by definition $c_j = c_0^{\;j}$, and the right side above is still an upper bound to $|c_j|^2$.  Hence everywhere in $\omega$,
\begin{align*}
W^{e}(\nabla' y_j |c_j, n_j ) &\leq c\left( |\nabla' y_j|^2 + |c_j|^2 + 1 \right) \nonumber \\
&\leq C\left(|\nabla' y_j|^2 + \frac{1}{|\adj \nabla' y_j|^2} + 1\right) \nonumber \\
&\leq C\left( W_{2D}(\nabla' y_j) + 1\right),
\end{align*}
\noindent using the growth in (\ref{eq:WeG}), the bound above and the coercivity in (\ref{eq:W2DG}).  This implies the bound
\begin{align}\label{eq:UB3}
\int_{\omega} W^e(\nabla ' y_j|c_j,n_j) dz' \leq \int_{\omega} W_{2D}(\nabla' y_j) dz' + E^2_j,
\end{align}
\noindent where recalling (\ref{eq:UB2}) and (\ref{eq:Nj}), the remainder $E^2_j$ is given by 
\begin{align}\label{eq:E2j}
E^2_j := C\left( \int_{\Gamma_{2\eta_j}} \left\{ W_{2D}( \nabla' y_j) + 1 \right\} dz' + \frac{1}{j}\right).
\end{align}

To recap, from (\ref{eq:UB1}) and (\ref{eq:UB3}), the entropic part of the energy is bounded above by the estimate 
\begin{align}\label{eq:Recap}
\int_{\Omega} W^e(\nabla_h y_j^{\;h}, n_j^{\;h}) \leq \int_{\omega} W_{2D}(\nabla' y_j) dz' + E^1_{h,j} + E^2_j.
\end{align}
\noindent It remains to bound the elasticity of the nematic mesogens.  

Consider the second term of $\mathcal{\tilde{I}}^h$ in (\ref{eq:UBIh}).  Our deformations and director fields have sufficient regularity, so 
\begin{align}\label{eq:boundIntegrand}
 \frac{\kappa_h}{h^2} \int_{\Omega} |(\nabla n_j^{\;h})(\adj \nabla y_j^{\;h})|^2 dz 
&= \kappa_h \int_{\Omega} |(\nabla' n_j|0)(\adj \nabla_h y_j^{\;h})|^2 dz.
\end{align} 
\noindent Here, we used the identity $(1/h)(\nabla n)(\adj \nabla y) = (\nabla_h n)(\nabla_h y)$ and the definition $n_j^{\;h}(z) := n_j(z')$.  We bound the integrand by a constant independent of $h$.  To do this, we first consider the pointwise estimate in (\ref{eq:bGrad}).  An application of the reverse triangle inequality on this bound yields for small $h$ the pointwise estimate
\begin{align}\label{eq:bound2}
| \partial_1 y_j^{\;h}(z)|^2 + |\partial_2 y_j^{\;h}(z)|^2 + \frac{1}{h^2}|\partial_3 y_j^{\;h}(z)|^2& = |\nabla_h y_j^{\;h}(z)|^2 \nonumber \\ 
&\leq \left( C_j h + |(\nabla' y_j|c_j)(z')|\right)^2  \nonumber \\
& \leq (\tilde{M}_j/3)^{1/2}, \;\;\;\; z \in \Omega .
\end{align}
\noindent Here $\tilde{M}_j$ is a constant which depends only on $y_j$ and $c_j$.  Then $F = (f_1|f_2|f_3) \in \mathbb{R}^{3\times3}$ satisfies
\begin{align*}
|\adj F|^2 = |\cof F|^2 &= |(f_2 \times f_3|f_3 \times f_1| f_1 \times f_2)|^2 \\
&= |f_2 \times f_3|^2 + |f_3 \times f_1|^2 + |f_1 \times f_2|^2 \\
&\leq |f_2|^2|f_3|^2 + |f_3|^2|f_1|^2 + |f_1|^2|f_2|^2,
\end{align*} 
\noindent and so we can bound from above (\ref{eq:boundIntegrand}),
\begin{align*}
\kappa_h \int_{\Omega} &|(\nabla' n_j|0)(\adj \nabla_h y)|^2 dz  \leq \kappa_h \int_{\Omega}| \nabla' n_j|^2 |\adj \nabla_h y_j^{\;h}|^2 dz \\&\leq  \kappa_h \int_{\Omega} | \nabla' n_j|^2 \left(\frac{1}{h^2}|\partial_2 y_j^{\;h}|^2 |\partial_3 y_j^{\;h}|^2 + \frac{1}{h^2} |\partial_3 y_j^{\;h} |^2 |\partial_1 y_j^{\;h}|^2 + |\partial_1 y_j|^2|\partial_2 y_j|^2 \right) dz  .
\end{align*}
\noindent  Applying the bound in (\ref{eq:bound2}) to this estimate, we conclude as desired 
\begin{align}\label{eq:UB5}
\frac{\kappa_h}{h^2} \int_{\Omega}|(\adj \nabla n_j^{\;h})(\nabla y_j^{\;h})|^2dz  \leq \kappa_h \tilde{M}_j \int_{\Omega} |\nabla' n_j|^2 dz =: \kappa_h M_j.
\end{align}
Here $M_j$ is a constant depending only on $y_j$, $c_j$ and $n_j$.

To complete the proof of Proposition \ref{1407251610}, it remains to show that in the limit as $h \rightarrow 0$, the energy is bounded as in (\ref{eq:PropEqn}).  From (\ref{eq:Recap}) and (\ref{eq:UB5}),
\begin{align}\label{eq:UB6}
\tilde{\mathcal{I}}^h(y_j^{\;h},n_j^{\;h}) \leq \int_{\omega} W_{2D}(\nabla' y_j) dz' + E^1_{h,j} + E^2_j + \kappa_h M_j.  
\end{align}
\noindent We now fix $j \in \mathbb{N}$ and take the limit as $h \rightarrow 0$.  Notice from (\ref{eq:yjhCon}), $\| \nabla_h y_j^{\;h} - (\nabla' y_j| c_j) \|_{L^2} \rightarrow 0$ as $h \rightarrow 0$.  This implies $\|\nabla_h y_j^{\;h} \|_{L^2} \leq C_j$ for some constant $C_j$ independent of $h$.  With these two observations, we conclude $E^1_{h,j} \rightarrow 0$ as $h \rightarrow 0$, see (\ref{eq:E1hj}).  Further, since $\kappa_h \rightarrow 0$ as $h \rightarrow 0$, $\kappa_h M_j \rightarrow 0$ since $M_j$ is independent of $h$.  Collecting these results and combining with (\ref{eq:UB6}), 
\begin{align*}
\limsup_{h \rightarrow 0} \tilde{\mathcal{I}}^h(y_j^{\;h},n_j^{\;h})   &\leq \limsup_{h \rightarrow 0}  \left(\int_{\omega} W_{2D}(\nabla' y_j) dz' + E^1_{h,j} + E^2_j + \kappa_h M_j \right) \nonumber \\
&= \int_{\omega} W_{2D}(\nabla' y_j) dz' + E^2_j.
\end{align*}
\noindent Finally, using (\ref{eq:seqWqc}), the fact that $| \Gamma_{2\eta_j}| \rightarrow 0$ as $j \rightarrow \infty$ ($\eta_j \rightarrow 0$, see Lemma \ref{1407251613} (ii), and (\ref{eq:E2j})) we conclude 
\begin{align*}
\limsup_{j \rightarrow \infty} \limsup_{h \rightarrow 0} \tilde{\mathcal{I}}^h(y_j^{\;h},n_j^{\;h})  &\leq \limsup_{j \rightarrow \infty} \left( \int_{\omega} W_{2D}(\nabla' y_j) dz' + E^2_j\right) \nonumber \\
& \leq \int_{\omega} W_{2D}^{qc}(\nabla' y) dz'.  
\end{align*}
\noindent We now choose a diagonal sequence $h_j \rightarrow 0$ as $j \rightarrow \infty$ so that this estimate is satisfied and $y^{h_j} \rightharpoonup y$ in $H^1(\Omega,\mathbb{R}^3)$.  This completes the proof.  
\endproof

\section{Effective energy of nematic elastomer membranes}\label{EffEn}

In this section, we provide an explicit formula for the effective energy of  
nematic elastomer membranes.

\subsection{Effective energy}

\begin{figure}
\centering%
\includegraphics[width=4in]{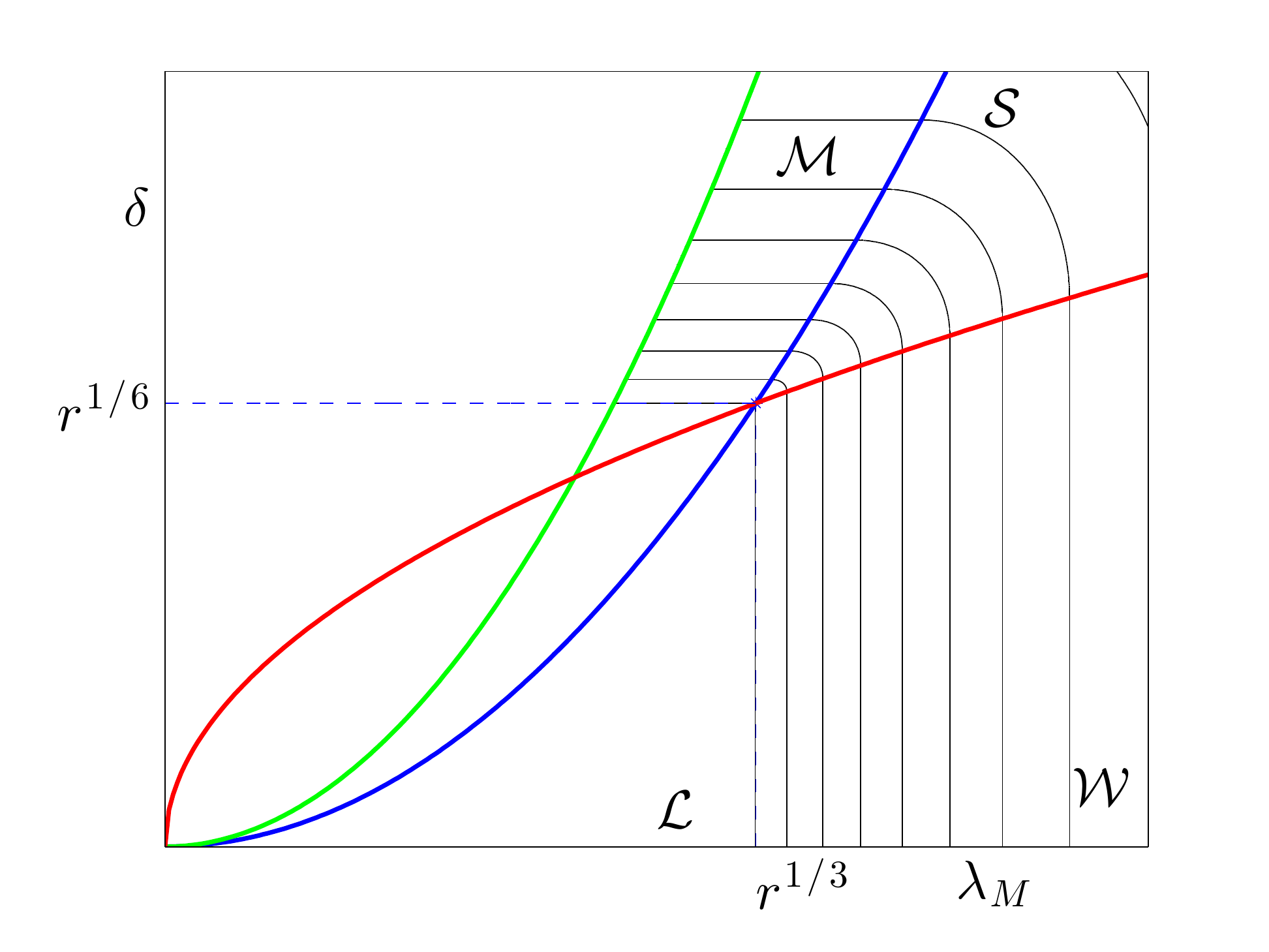}
\caption{The effective energy of nematic elastomer membranes.  $\lambda_M$ is the largest principal stretch and 
$\delta$  is the areal stretch.  The energy is zero in the region marked $\mathcal L$.}\label{fig:eff}
\end{figure}

\begin{theorem}\label{1407231520}
Let $W^e$ as in (\ref{eq:We}).  For any $\tilde F \in {\mathbb R}^{3\times 2}$, let 
\begin{eqnarray} \label{eq:W2D}
W_{2D} (\tilde{F}) = \inf_{c \in {\mathbb R}^3, n \in \mathbb{S}^2} W^e ( (\tilde{F}|c), n).
\end{eqnarray}
Then, the effective energy of the nematic membrane $W_{2D}^{qc}  \equiv W^{mem}$ is given by
\begin{eqnarray} \label{eq:Wmem}
W^{mem}(\tilde{F}) =\frac{\mu}{2} \begin{dcases}
0 &\text{ if } (\lambda_M, \delta) \in \mathcal{L}, \\
r^{1/3} \left(\frac{2\delta(\tilde{F})}{r^{1/2}} + \frac{1}{\delta(\tilde{F})^2} \right)- 3 &\text{ if } (\lambda_M, \delta) \in \mathcal{M}, \\
r^{1/3} \left(\frac{\lambda_M(\tilde{F})^2}{r} + \frac{2}{\lambda_M(\tilde{F})}\right) -3 &\text{ if } (\lambda_M,\delta) \in \mathcal{W}, \\
r^{1/3} \left(\frac{\lambda_M(\tilde{F})^2}{r} + \frac{\delta(\tilde{F})^2}{\lambda_M(\tilde{F})^2} + \frac{1}{\delta(\tilde{F})^2} \right) - 3 &\text{ if } (\lambda_M,\delta) \in \mathcal{S}. 
\end{dcases}
\end{eqnarray}
Here, $0 \le \lambda_m(\tilde F) \le \lambda_M(\tilde F)$ are the singular values of $\tilde F$ (i.e., the eigenvalues of $(\tilde F^T \tilde F)^{1/2}$), $\delta(\tilde F) = \lambda_m(\tilde F)\lambda_M(\tilde F) = \sqrt{\det \tilde F^T \tilde F}$ and 
\begin{eqnarray}
\mathcal{L} &:=& \{(\lambda_M, \delta) : \lambda_M^2 \geq \delta, \lambda_M \leq r^{1/3}, \delta \leq r^{1/6} \}, \\
\mathcal{W} &:=& \{ (\lambda_M, \delta) : \lambda_M > r^{1/3}, \delta < \lambda_M^{1/2} \}, \\
\mathcal{M} &:=& \{ ( \lambda_M, \delta) : \delta > r^{1/6} , r^{-1/2}\lambda_M^2 < \delta \leq \lambda_M^2 \},  \\
\mathcal{S} &:=& \{(\lambda_M, \delta) : \lambda_M^{1/2} \leq \delta \leq r^{-1/2} \lambda_M^2 \} .
\end{eqnarray}
\end{theorem}

\begin{remark}
Some care needs to be taken when dealing with extended real-valued quasiconvex functions. 
Indeed, the fact that a function $f:\mathbb{R}^{3\times 2}\to\mathbb{R}\cup\{+\infty\}$ is quasiconvex (according to the defintion of Section \ref{1407231507} of this paper) does not imply that the associated functional $\int_{\omega}f(\nabla 'y')dx'$ is  sequentially weak$^*$ lower semicontinuous on $W^{1,\infty}(\omega,\mathbb{R}^3)$ \cite{BallMurat}. In the current situation, thanks to Remark \ref{1407231441}, the relaxed energy density has polynomial growth and therefore weak lower semincontinuity is true for the relaxed functional.
Alternatively, we refer the interested readers to Ball and James \cite{BaJa15} where a more restrictive definition of quasiconvexity for extended real value functions is presented. This definition guarantees weak lower semicontinuity of functionals associated to extended real value integrand functions.
It is an easy computation to show that both the approach pursued in what follows and the  relaxation technique based on the alternative definition of quasiconvex envelope  give the same results for the functionals considered in this paper.

\end{remark}

\noindent {\it Proof of Theorem \ref{1407231520}.}
Recall that the quasiconvex envelope of an extended value function is not in general bounded from above by the rank-one convex envelope.  However, we show that this bound is true for $W_{2D}$.  By Remark \ref{1407231441},  $W_{2D}^{qc}$ is a finite-valued, a quasiconvex function and $W_{2D}^{qc}=(W_{2D}^{qc})^{qc}$. Therefore, if we substitute $f=W_{2D}^{qc}$ in $(\ref{1303082354})$ we obtain
\begin{eqnarray}\label{1303082355}
\bigl(W_{2D}^{qc}\bigr)^{pc}\leq \bigl(W_{2D}^{qc}\bigr)^{qc}\leq \bigl(W_{2D}^{qc}\bigr)^{rc}.
\end{eqnarray}
Then, by $(\ref{1303082357})$ we conclude
\begin{eqnarray}
W_{2D}^{pc}\equiv \bigl(W_{2D}^{pc}\bigr)^{pc}\leq\bigl(W_{2D}^{qc}\bigr)^{pc}\leq \bigl(W_{2D}^{qc}\bigr)^{qc}\leq \bigl(W_{2D}^{qc}\bigr)^{rc}\leq W_{2D}^{rc}
\end{eqnarray}
and recover the classical inequality 
\begin{eqnarray}\label{1309182125}
(W_{2D})^{pc}\leq (W_{2D})^{qc}\leq (W_{2D})^{rc}.
\end{eqnarray}

We show in Lemma \ref{lem:W2D} and Lemma \ref{lem:ub} that $W^{mem} \le W_{2D}^{pc}$.  We show in 
Lemma \ref{lem:lb} that $W_{2D}^{rc} \le W^{mem}$.    Combining these with (\ref{1309182125}), 
\begin{eqnarray}\label{1407231457}
W^{mem} \le W_{2D}^{pc} \le W_{2D}^{qc} \le W_{2D}^{rc} \le W^{mem},
\end{eqnarray}
and the result follows.
\endproof

\begin{figure}
        \centering
        \begin{subfigure}[b]{0.4\textwidth}
\includegraphics[width=2.75in]{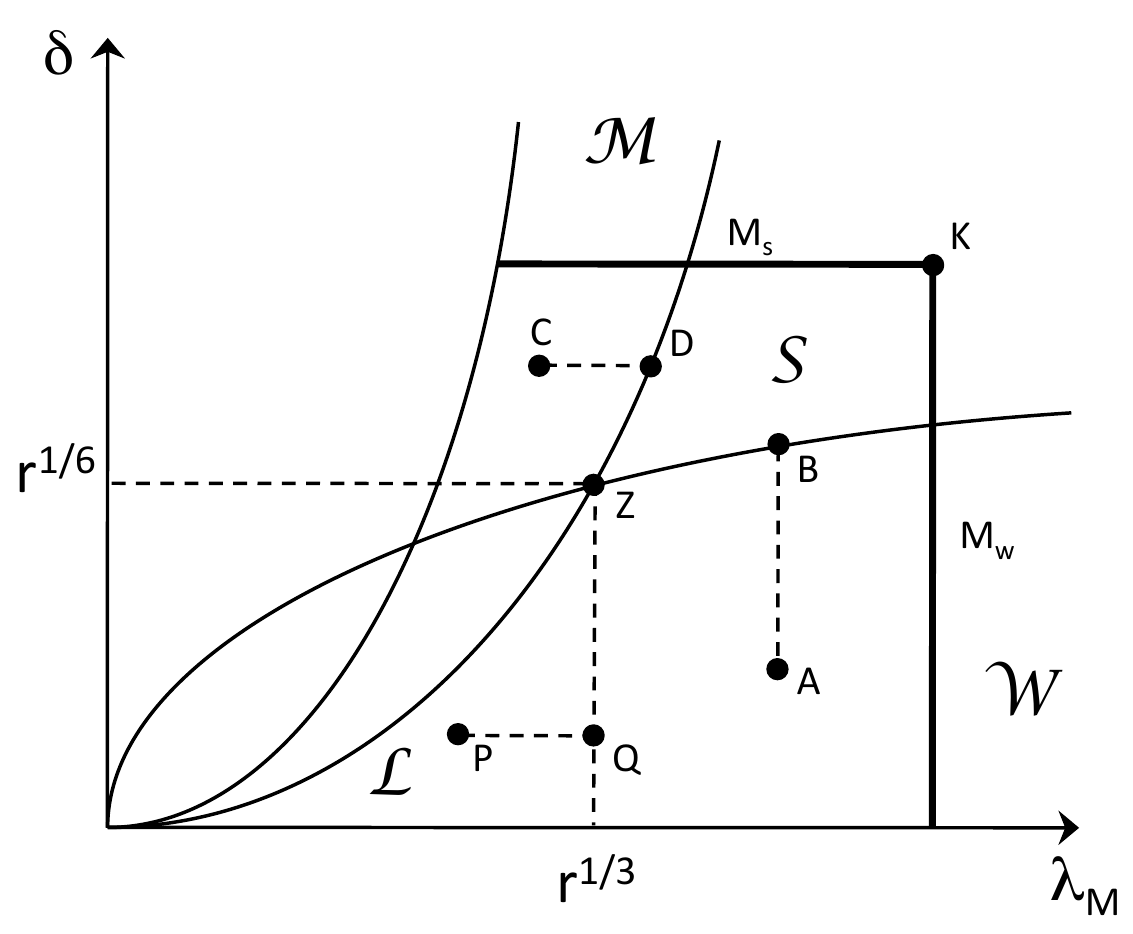}\quad
                \caption{}
                \label{fig:a}
        \end{subfigure}%
\qquad
        ~ 
        \begin{subfigure}[b]{0.4\textwidth}
\includegraphics[width=7.5cm]{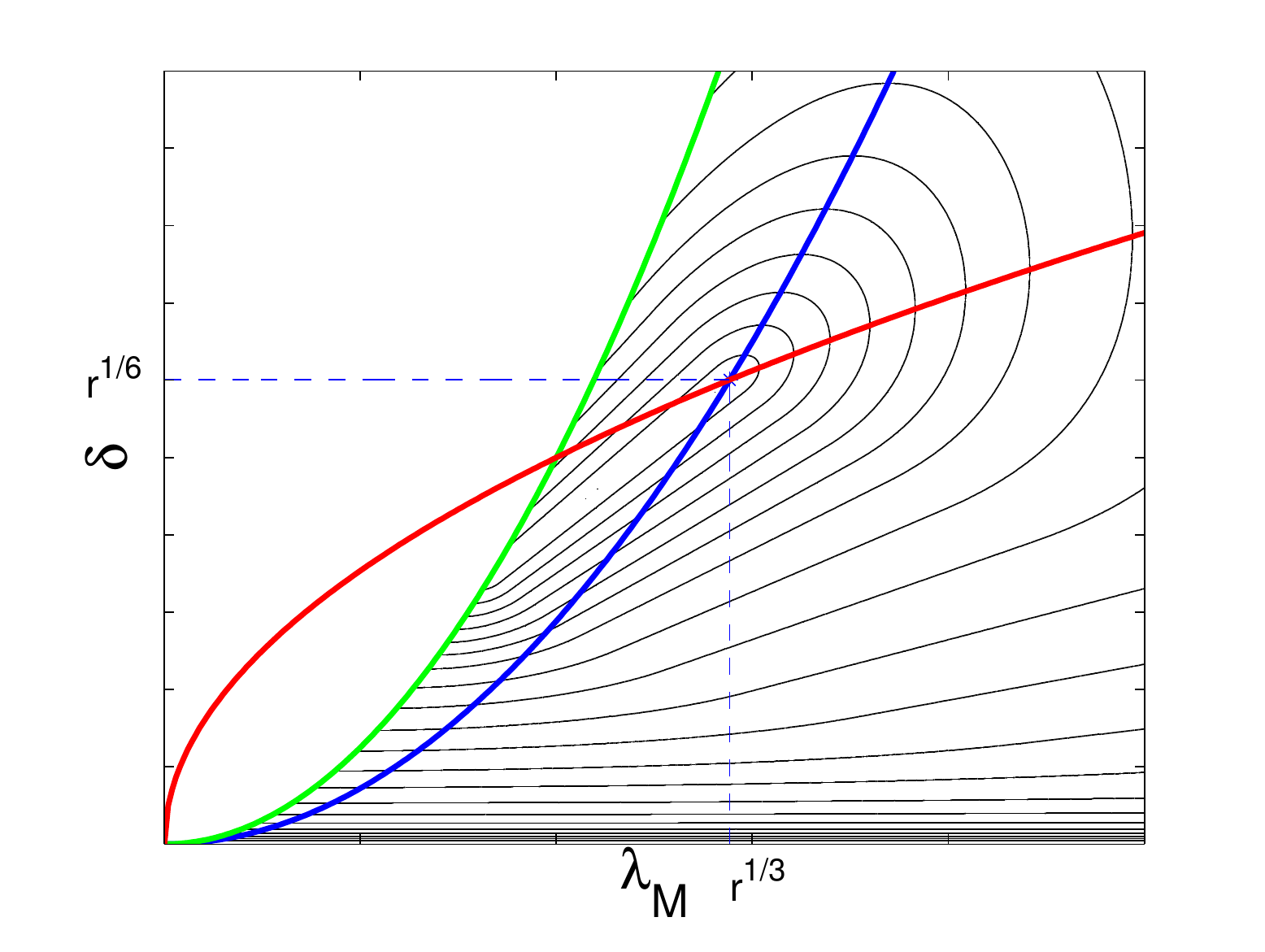}%
                \caption{}
                \label{fig:b}
        \end{subfigure}
\caption{
$(a)$: Idea of the proof of Lemma \ref{lem:lb}.
$(b)$: level curves of $W_{2D}$ in the space $(\lambda_M,\delta)$.
The web version of this article contains the above plot figures in color.}
\label{fig:lam}
\end{figure}

\subsection{Step 1: A formula for $W_{2D}$} \label{sec:W2D} 

\begin{lemma} \label{lem:W2D}
For $W_{2D}$ defined in (\ref{eq:W2D}),
\begin{eqnarray}  \label{eq:W2Dmin}
W_{2D} (\tilde F) = \begin{dcases}
\min_{i \in \{1, \dots, 3\}} \varphi_i (\lambda_M (\tilde F), \delta (\tilde F)), &\text{ if } \rank \tilde{F} = 2 \\
+\infty &\text{ otherwise }
\end{dcases}
\end{eqnarray}
where

\begin{eqnarray}\label{1407221344}
\min_{i \in \{1, \dots, 3\}} \varphi_i (\lambda_M , \delta ) = \begin{dcases}
\varphi_1 (\lambda_M , \delta) &\text{ if }  \lambda_M\delta\geq r^{1/2}, \delta\leq\lambda_M^2\\
\varphi_2 (\lambda_M , \delta) &\text{ if }  \lambda_M\delta\leq r^{-1/2},\delta\leq\lambda_M^2\\
\varphi_3 (\lambda_M, \delta) &\text{ if }  \lambda_M\delta\in (r^{-1/2},r^{1/2}),\delta\leq\lambda_M^2
\end{dcases}
\end{eqnarray}

with
\begin{eqnarray}\label{eq:varphi1}
\varphi_1(\lambda_M,\delta)&:=&\displaystyle{
\frac{\mu}{2}\Bigl\{r^{1/3} \Bigl[\frac{\lambda_M^2}{r}+\frac{\delta^2}{\lambda_M^2}  +\frac{1}{\delta^2} 
\Bigr]-3\Bigr\}
},\\
\label{eq:varphi2}
\varphi_2(\lambda_M,\delta) &:=&
\displaystyle{
\frac{\mu}{2}\Bigl\{r^{1/3} \Bigl[\lambda_M^2+\frac{\delta^2}{\lambda_M^2}+\frac{1}{r\delta^2} 
\Bigr]-3\Bigr\}
},\\
\label{eq:varphi3}
\varphi_3(\lambda_M,\delta)&:=&\displaystyle{\begin{dcases}
\frac{\mu}{2}\Bigl\{r^{1/3} \Bigl[\frac{\delta^2}{\lambda_M^2}+2\frac{\lambda_M}{r^{1/2}\delta} 
\Bigr]-3\Bigr\}  &\text{ if } \lambda_M \delta \in (r^{-1/2},r^{1/2}), \\
+\infty &\text{ otherwise }.
\end{dcases}} 
\end{eqnarray}

\end{lemma}

\noindent {\it Proof.}
The proof is an explicit calculation.  To begin, if $\rank \tilde{F} \neq 2$, then $\det(\tilde{F}|c) = 0$ for every $c \in \mathbb{R}^3$.  This implies $W^e(\tilde{F}|c,n) = +\infty$ for every $c \in \mathbb{R}^3$. Then $W_{2D}(\tilde{F}) = +\infty$.  Thus for the remainder of this section, we restrict our attention to the case that $\rank \tilde{F} = 2$.  

Let $f_{\tilde{F},n}(c) := W_0(\tilde{F}|c,n)$ and $g_{\tilde{F}}(c) := c^T \adj \tilde{F} - 1$.  Then, $\inf_{c \in \mathbb{R}^3} W^e(\tilde{F}|c,n)$ is equivalent to the optimization
\begin{align*}
\inf_{c \in \mathbb{R}^3} \left\{ f_{\tilde{F},n}(c) : g_{\tilde{F}}(c) = 0 \right\}.
\end{align*}
\noindent Here $f_{\tilde{F},n}$ is a convex, differentiable function and $g_{\tilde{F}}$ is an affine equality constraint.   It follows that $c_0$ is a global minimizer of this optimization if and only if there exists a $\lambda \in \mathbb{R}$ such that $\nabla f_{\tilde{F},n}(c_0) + \lambda \nabla g_{\tilde{F}}(c_0) = 0$ (see for instance \cite{BV04}, Section 5.5.3).  Solving this equation, we obtain 
\begin{align*}
c_0 = \frac{\ell \adj \tilde{F}}{|\ell^{1/2} \adj \tilde{F}|^2}
\end{align*}
\noindent for $\ell$ defined in (\ref{eq:stepLength}).  

Let $\widetilde{W}(\tilde{F},n) := W_0(\tilde{F}|c_0,n)$. Since $c_0$ is a global minimizer for the constrained optimization above, $\widetilde{W}(\tilde{F},n) = \inf_{c \in \mathbb{R}^3} W^e(\tilde{F}|c,n)$. Then from (\ref{eq:W2D}), it follows that $\inf_{n \in \mathbb{S}^2} \widetilde{W}(\tilde{F},n) = W_{2D}(\tilde{F})$.  For this optimization, we simplify the analysis through a change of variables.  We write $\tilde{F} = QDR$ for $Q \in SO(3)$, $R \in O(2)$ and $D$ a diagonal matrix as in (\ref{eq:pd}) with $\lambda_M \geq \lambda_m > 0$ as the singular values.  We can say $\lambda_m > 0$ since $\rank \tilde{F} = 2$.  Additionally, we set $n = Qm$, and impose the $\mathbb{S}^2$ constraint via $m_3^2 = 1 - m_1^2 - m_2^2$.  Then by direct substitution,
\begin{align*}
\widetilde{W}(QDR,Qm) = \frac{\mu}{2} \left\{ r^{1/3}\left[\gamma \lambda_M^2 - \xi_2 (\lambda_M^2-\lambda_m^2) + \frac{1}{r(\gamma -1)\lambda_M^2\lambda_m^2}\right]-3\right\} =: \widetilde{\varphi}(\lambda_M, \lambda_m, \gamma,\xi_2),
\end{align*}
\noindent where
\begin{align*} 
\gamma = \xi_1 + \xi_2, \;\;\;\; \xi_i(m_i) = 1-\alpha m_i^2, \;\;\;\;  i = 1,2, \;\;\;\; \alpha = \frac{r-1}{r}.
\end{align*}
\noindent Here $\alpha \in [0,1)$ since $r \geq 1$.  Further, we let $\delta = \lambda_M \lambda_m$, and set 
\begin{align}\label{eq:varphi}
\varphi(\lambda_M,\delta,\gamma,\xi_2) &:= \widetilde{\varphi}(\lambda_M,\delta/\lambda_M,\gamma,\xi_2) \nonumber \\
&\;=  \frac{\mu}{2} \left\{ r^{1/3}\left[\gamma \lambda_M^2 - \xi_2\left( \lambda_M^2 - \left(\frac{\delta}{\lambda_M}\right)^2\right) + \frac{1}{r(\gamma -1)\delta^2}\right]-3\right\}. 
\end{align}
\noindent Note that the constraint $\lambda_M \geq \lambda_m> 0$ implies $\lambda_M^2 \geq \delta > 0$.  

$\widetilde{W}$ is dependent on only four constrained variables.  Consider the closed set 
\begin{align*}
\mathcal{B} := \left\{ (\gamma,\xi_2) : \; \gamma - \xi_2 \leq 1,\; \xi_2 \in [1-\alpha,1] \text{ and } \gamma \in [2-\alpha,2] \right\}.
\end{align*}
\noindent $\mathcal{B}$ combined with the constraint $\lambda_M^2 \geq \delta > 0$ give the admissible set for $\varphi$.  Hence, we have 
\begin{align*}
W_{2D}(\tilde{F}) = \inf_{n \in \mathbb{S}^2} \widetilde{W}(\tilde{F},n) = \inf_{\gamma, \xi_2} \left\{ \varphi(\lambda_M(\tilde{F}), \delta(\tilde{F}),\gamma, \xi_2) : \lambda_M^2 \geq \delta > 0, (\gamma, \xi_2) \in \mathcal{B} \right\}. 
\end{align*}
\noindent Observe that 
\begin{align*}
\inf_{\xi_2}\left\{ \varphi(\lambda_M,\delta, \gamma, \xi_2) : \lambda_M^2 \geq \delta > 0 , (\gamma,\xi_2)\in \mathcal{B} \right\} = \varphi(\lambda_M, \delta, \gamma,1) =: \varphi_0(\lambda_M,\delta,\gamma)
\end{align*}
\noindent by (\ref{eq:varphi}) since the second term in the brackets is non-positive.  Then,
\begin{align*}
W_{2D}(\tilde{F}) = \inf_{\gamma} \left\{ \varphi_0(\lambda_M(\tilde{F}), \delta(\tilde{F}),\gamma ): \lambda_M^2 \geq \delta > 0, \gamma \in [2-\alpha,2] \right\},
\end{align*}
\noindent where 
\begin{align*}
\varphi_0(\lambda_M, \delta, \gamma) = \frac{\mu}{2} \left\{ r^{1/3}\left[(\gamma - 1) \lambda_M^2 + \left(\frac{\delta}{\lambda_M}\right)^2+ \frac{1}{r(\gamma -1)\delta^2}\right]-3\right\}.
\end{align*}

$\varphi_0$ is a continuous function on this constrained set (which is moreover bounded in $\gamma$).  It is also differentiable for $\gamma$ in the open domain $(2-\alpha,2)$.  It follows that the infimum is attained.  Further, $\bar{\gamma}$ minimizes $\varphi_0$ only if it is on the boundary, i.e $\bar{\gamma} = 2-\alpha$ or $\bar{\gamma} = 2$, or it is a critical point, i.e. $\partial_{\gamma} \varphi_0(\bar{\gamma}) = 0$ and $ \bar{\gamma} \in (2-\alpha,2)$.  

We proceed case by case.  Letting $\bar{\gamma} = 2- \alpha$, observe $\varphi_0(\lambda_M,\delta,2-\alpha) = \varphi_1(\lambda_M,\delta)$ in (\ref{eq:varphi1}).  For the other boundary $\bar{\gamma} = 2$, we obtain $\varphi_0(\lambda_M,\delta,2) = \varphi_2(\lambda_M,\delta)$ in (\ref{eq:varphi2}).  Finally, in computing the critical point $\partial_{\gamma} \varphi_0(\bar{\gamma}) = 0$, we obtain 
\begin{align*}
\bar{\gamma} = \frac{r^{-1/2}}{\lambda_M \delta} + 1.
\end{align*}
\noindent Direct substitution $\varphi_0(\lambda_M, \delta,\bar{\gamma})$ yields the equation given for the finite portion of $\varphi_3$ in (\ref{eq:varphi3}).  Recalling that $\bar{\gamma}$ must lie in the domain $(2-\alpha,2)$, we set $\varphi_3 = +\infty$ if $\bar{\gamma}$ does not lie in this set.  This is the full result in (\ref{eq:W2Dmin}).  

To complete the proof we compute the minimum in Eq. (\ref{eq:W2Dmin}) thus yielding
(\ref{1407221344}). First, observe that
\begin{eqnarray*}
\varphi_1-\varphi_3=\frac{\mu}{2}r^{1/3}\Bigl(r^{-1/2}\lambda_M-\frac{1}{\delta}\Bigr)^2\geq 0\Longrightarrow\varphi_3\leq\varphi_1,\\
\varphi_2-\varphi_3=\frac{\mu}{2}r^{1/3}\Bigl(\lambda_M-\frac{r^{-1/2}}{\delta}\Bigr)^2\geq 0\Longrightarrow\varphi_3\leq\varphi_2,
\end{eqnarray*}
thus proving (\ref{1407221344}) in the region $\lambda_M\delta\in(r^{-1/2},r^{1/2})$. To complete the computation in the remaining regions observe that
$$
\varphi_1-\varphi_2=\frac{\mu}{2}r^{1/3}\Bigl(1-\frac{1}{r}\Bigr)\Bigl(\frac{1}{\delta}+\lambda_M\Bigr)\Bigl(\frac{1}{\delta}-\lambda_M\Bigr),
$$
yielding $\varphi_1\leq\varphi_2$ if $\delta\geq\lambda_M^{-1}$ and 
$\varphi_2\leq\varphi_1$ if $\delta\leq\lambda_M^{-1}$. Therefore we have
$$
\min_{1,\dots,3}\varphi_i(\lambda_M,\delta)=\varphi_1 \text{ if } \lambda_M\delta\geq r^{1/2},\delta\leq\lambda_M^2
$$
and
$$
\min_{1,\dots,3}\varphi_i(\lambda_M,\delta)=\varphi_2  \text{ if  }  \lambda_M\delta\leq r^{-1/2},\delta\leq\lambda_M^2
$$
as required.
\endproof

\subsection{Step 2: Upper bound or $W^{mem} \le W_{2D}^{pc}$} \label{sec:ub}

\begin{lemma} \label{lem:ub}
Let $W^{mem}$ be as in (\ref{eq:Wmem}) and $W_{2D}$ as in (\ref{eq:W2D}).  Then,
for each $\tilde F \in {\mathbb R}^{3\times 2}$,
\begin{eqnarray}
W^{mem} (\tilde F) \le W_{2D}^{pc} (\tilde F).
\end{eqnarray}
\end{lemma}

\noindent {\it Proof.}  We prove this in two parts.  In Part 1, we prove that $W^{mem}$ is polyconvex and in Part 2 we prove that $W^{mem} \le W_{2D}$.  The result follows.

\noindent Part 1.  We now show that $W^{mem}$ is polyconvex.  
First, observe from (\ref{eq:Wmem}) that there exists a function $\psi: {\mathbb R}_+^2 \rightarrow {\mathbb R}$ (here by $\mathbb{R}_+$ we denote the set of all non-negative real numbers) such that
\begin{eqnarray} \label{eq:psi}
W^{mem}(\tilde F) = \psi(\lambda_M (\tilde F), \delta (\tilde F)).
\end{eqnarray}
It also follows by verification (also see  Proposition 2 of DeSimone and Dolzmann \cite{DD00}) that $\psi$ is convex and $\psi$ is non-decreasing in each argument (i.e., $\psi(s,t)$ in nondecreasing in $s$ for fixed $t$ and nondecreasing in $t$ for fixed $s$.  We then notice that $\lambda_M (\tilde{F}) = \sup_{m \in \mathbb{S}^1} |\tilde{F}m|$ is convex in $\tilde{F}$.  Further, $\delta (\tilde F) = | \adj \tilde F| $ is convex in $\adj \tilde F$.  Since the composition of convex function with a non-decreasing and convex function results in a convex function, we conclude that there exists a  convex function $g: {\mathbb R}^{3\times2} \times {\mathbb R}^3 \rightarrow \mathbb{R} $ such that
$$
\psi(\lambda_M (\tilde F), \delta (\tilde F)) = g(\tilde F, \adj \tilde F).
$$
Combining with (\ref{eq:psi}), we conclude that
\begin{eqnarray}
W^{mem}(\tilde F) =  g(\tilde F, \adj \tilde F)
\end{eqnarray}
for convex $g$.  By definition of polyconvexity, $W^{mem}$ is polyconvex.

\noindent Part 2.  We now show that $W^{mem} \le W_{2D}$.  We show by explicit calculation in the Appendix A that
$$
W^{mem}(\tilde F) \le \varphi_i (\lambda_m (\tilde F), \delta (\tilde F)), \quad i = 1, \dots 3.
$$
It follows that 
$$
W^{mem}(\tilde F) \le \min_{i \in \{1, \dots, 3\}} \varphi_i (\lambda_m (\tilde F), \delta (\tilde F)) = W_{2D} (\tilde F)
$$
from (\ref{eq:W2Dmin}).
\endproof

\subsection{Step 3: Lower bound or $W_{2D}^{rc} \le W^{mem}$} \label{sec:lb}

\begin{lemma} \label{lem:lb}
Let $W^{mem}$ be as in (\ref{eq:Wmem}) and $W_{2D}$ as in (\ref{eq:W2D}).  Then,
for each $\tilde F \in {\mathbb R}^{3\times 2}$,
\begin{eqnarray}
W_{2D}^{rc} (\tilde F) \le W^{mem} (\tilde F).
\end{eqnarray}
\end{lemma}
The proof makes repeated use of lamination.  We collect the calculations in the following proposition.
\begin{proposition} \label{prop:lam}
Let $q,d \in {\mathbb R}$ with $q >0$, $q^2\ge d$ and define
\begin{eqnarray*}
\calK&:=& \bigl\{\tilde G \in{\mathbb R}^{3\times 2}: \lambda_{M}(\tilde G)=q,\,\delta(\tilde G)=d \bigr\}, \\
M_w &:= & \bigl\{ \tilde G \in{\mathbb R}^{3\times 2}: \lambda_{M}(\tilde G) = q, \, \delta(\tilde G) \in [0,d] \bigr\}, \\
M_s& := & \bigl\{\tilde G \in{\mathbb R}^{3\times 2}: \lambda_{M}(\tilde G) \in [d^{1/2},q], \,\delta(\tilde G)= d \bigr\} .
\end{eqnarray*}
Then, for $d>0$ $M_w \subset \calK^{(1)}$ and for $d\geq 0$ $M_s\subset \calK^{(1)}$.
\end{proposition}
\noindent {\it Proof.} We begin with $M_w$.  Let $ \tilde G \in M_w$ with $\lambda_{M}(\tilde G) = q, \, \delta(\tilde G) =\bar{\delta} \in [0,d]$.  Using the polar decomposition theorem, we can take
\begin{equation*} \label{eq:gtilde}
 \tilde G = \left(
\begin{array}{ccc}
q & 0\\
0 & \frac{\bar{\delta}}{q}  \\
0 & 0 
\end{array} \right).
\end{equation*}
Define 
$$
 \tilde G^\pm = \left(
\begin{array}{ccc}
q & 0   \\
0 & \pm \frac{d}{q} \\
0 & 0
\end{array} \right); \quad
\theta = \frac{1}{2} \left(1+\frac{\bar{\delta}}{d} \right) .
$$
Note that $ \tilde G^\pm \in \calK$, $\theta \in [0,1]$ since $\bar{\delta} \le d$,  rank $( \tilde G^+ - \tilde G^-) = 1$ and
$ \tilde G = \theta \tilde G^+ + (1-\theta) \tilde G^-$.  Therefore, $M_w \subset \calK^{(1)}$.
 
 The proof of $M_s \subset \calK^{(1)}$ is similar (also see \cite[Theorem 3.1]{DD00}). Again, using the polar decomposition theorem, we can take $\tilde G \in M_s$ as a diagonal matrix.   
First, let us assume $c\neq 0$ and $\sqrt{d} \leq c\leq q$ which corresponds to $c   \geq d/c$ and define  
\begin{align}\label{eq:gtildepm}
 \tilde G^\pm := \left(
\begin{array}{ccc}
c & \pm \xi   \\
0 &  \frac{d}{c} \\
0 & 0
\end{array} \right).
\end{align}
Note $\delta(\tilde G^\pm) = d$.  Further, the eigenvalues of $(\tilde G^\pm)^T  \tilde G^\pm$ are
$$
\frac{1}{2} \left(\xi^2+\frac{d^2}{c^2}+c^2\right) \pm \sqrt{\frac{1}{4}\left(\xi^2+\frac{d^2}{c^2}+c^2\right)^2-d^2}.
$$
So the choice 
$$
\xi^2
= \frac{d^2}{q^2}+q^2-\frac{d^2}{c^2}-c^2
= \frac{1}{q^2}\Bigl[
q^2-\frac{d^2}{c^2}\Bigr]\Bigl[q^2-c^2\Bigl]\geq 0
$$
makes $\lambda_M(\tilde G^\pm) =q$.  Therefore, $\tilde G^\pm \in \calK$.  Further,
$\tilde G = \frac{1}{2} \tilde G^+ + \frac{1}{2} \tilde G^-$ and rank $( \tilde G^+ - \tilde G^-) \le 1$.   For the case $\tilde{G} \in M_s$ such that $\tilde{G} = 0$, replace the diagonal entries in (\ref{eq:gtildepm}) with $0$ and repeat the argument. Therefore, $M_s \subset \calK^{(1)}$.
\endproof
\vspace{\baselineskip}

\noindent {\it Proof of Lemma \ref{lem:lb}.} 
We show that
\begin{equation} \label{eq:wubrc}
W_{2D}^{rc}(\tilde F) \le W^{mem} (\tilde F)
\end{equation}
region by region.  In the region $\mathcal{S}$, note $W_{2D} = W^{mem}$ and the result follows.  

Now, let $\tilde F \in \mathcal{W}$ with  $q = \lambda_M(\tilde F )\ge r^{1/3}$ and $d=\delta(\tilde F)\leq q^{1/2}$.   This corresponds to the point $A$ in Figure \ref{fig:lam}.   Let 
$$
\tilde{\calK}=\bigl\{\tilde G \in{\mathbb R}^{3\times 2}: \lambda_{M}(\tilde G)=q,\,\delta(\tilde G)=q^{1/2}\bigr\}.
$$
This set corresponds to the point $B$ in Figure \ref{fig:lam}a.  By Proposition \ref{prop:lam}, we have 
$\tilde F \in \tilde{\calK}^{(1)}$.  Therefore,  there exists $\lambda \in [0,1]$ and $\tilde G_1, \tilde G_2 \in \tilde \calK$ with $\rank (\tilde{G}_1 - \tilde{G}_2) \leq 1$ such that 
\begin{eqnarray*}
W^{rc}_{2D}(\tilde F) &\le& \lambda W^{rc}_{2D} (\tilde G_1) + (1-\lambda) W^{rc}_{2D} (\tilde G_2) \\
&\le& \lambda W_{2D} (\tilde G_1) + (1-\lambda) W_{2D} (\tilde G_2)\\
&=& \frac{\mu}{2}\Bigl[r^{1/3}\Bigl(\frac{q^2}{r}+\frac{2}{q}\Bigr)-3\Bigr] \\
&=& W^{mem}(\tilde F).
\end{eqnarray*}
Above, the first two inequalities follow from the fact $W^{rc}_{2D}$ is rank-one convex and $W^{rc}_{2D}\le W_{2D}$. The following two equalities is by explicit verification of the formula.

Now, let $\tilde F \in \mathcal{M}$ with  $d=\delta(\tilde F) \ge r^{1/6}$ and $q=\lambda_M(\tilde F )\in[d^{1/2}, r^{1/4}d^{1/2}]$.  This corresponds to the point $C$ in Figure \ref{fig:lam}a.  Let
$$
\tilde{\calK}=\bigl\{\tilde G \in{\mathbb R}^{3\times 2}: \lambda_{M}(\tilde G)=r^{1/4}d^{1/2},\,\delta(\tilde G)=d \bigr\}.
$$
This set corresponds to the point $D$ in Figure \ref{fig:lam}a. Therefore, by Proposition \ref{prop:lam}, we have 
$\tilde F \in \tilde{\calK}^{(1)}$.  Therefore, arguing as before, $W^{rc}_{2D}(\tilde F) \le W^{mem}(\tilde F)$.

Finally, let $\tilde F \in \mathcal{L}$ with $q=\lambda_M(\tilde F )\le r^{1/3}$ and $d=\delta(\tilde F) \le \min\{q^2, r^{1/6}\}$.  This corresponds to the point $P$ in Figure \ref{fig:lam}a.  Let
$$
\tilde{\calK}=\bigl\{\tilde G \in{\mathbb R}^{3\times 2}: \lambda_{M}(\tilde G)=r^{1/3},\,\delta(\tilde G)=d \bigr\}
$$
and
$$
\calK=\bigl\{\tilde G \in{\mathbb R}^{3\times 2}: \lambda_{M}(\tilde G)=r^{1/3},\,\delta(\tilde G)=r^{1/6} \bigr\}.
$$
These sets correspond to the points $Q$ and $Z$ in Figure \ref{fig:lam}-(a) respectively.  From Proposition \ref{prop:lam}, 
$\tilde F \in \tilde{\calK}^{(1)}$.  Further, again by Proposition \ref{prop:lam}, $\tilde \calK \subset \calK^{(1)}$.  In other words, $\tilde F \in \calK^{(2)}$.  We can again argue as above to show that $W^{rc}_{2D}(\tilde F) \le 0 = W^{mem}(\tilde F)$ as required.

\endproof


\section{Characterization of fine-scale features}\label{Micro}

The energy density $W_{2D}$ is not quasiconvex.   Thus a membrane with this energy density is able to relax its energy to that of $W_{2D}^{qc}$ through the introduction of fine-scale features.  In this section, we characterize these features.   Briefly, we show that the features in region $\mathcal{M}$ are essentially planar involving oscillations of the director (i.e., no wrinkling) while those in $\mathcal{W}$ are necessarily wrinkles (i.e., uniform director).   Further, we show that there are no fine-scale features in region $\mathcal{S}$.

To characterize the fine-scale features, we consider the two-dimensional energy
\begin{align}\label{eq:I2D}
I_{2D}(y) = \int_{\omega} W_{2D}(\nabla' y) dx'
\end{align}
subject to affine boundary conditions, i.e the space of deformations $\mathcal{A}_{\bar{F}} := \{ y \in H^1(\omega,\mathbb{R}^3) \colon y - \bar{F}x' \in H^1_0(\omega,\mathbb{R}^3)\}$ with $\bar{F} \in \mathbb{R}^{3\times2}$.
It is known (cf. Lemma 3.1(ii) and Lemma 6.2, \cite{CD06}) that there exists weakly converging minimizing sequences  that satisfy
\begin{align}\label{eq:minSeq}
y_j \rightharpoonup \bar{F} x' \;\; \text{in} \;\; H^1(\omega,\mathbb{R}^3) \;\; \text{ with } \;\;  I_{2D}(y_j) \rightarrow \inf_{\mathcal{A}_{\bar{F}}} I_{2D} = |\omega| W^{mem}(\bar{F}) \:\: \text{as } \;\; j \rightarrow \infty.
\end{align}
Let $\nu_x$ be any  $H^1$ gradient Young measure generated by such a sequence.  Since $\{ y_j\}$ is a minimizing sequence for $I_{2D}$, it is also a minimizing sequence for the relaxation $\int_{\omega} W^{mem}(\nabla' y) dx'$.  Further, since $W^{mem}$ is non-negative and bounded as in (\ref{eq:W2DQCBound}), it follows from Theorem 1.3 of Kinderlehrer and Pedregal \cite{KP92} that
\begin{align}\label{eq:YMLimit}
f(\nabla y_j ) \rightharpoonup \langle \nu_x, f \rangle \;\; \text{ in } L^1(A) \;\; \text{ for any } f \in \left\{ g \in C(\mathbb{R}^{3\times2}) \colon \sup_{\tilde{G}\in \mathbb{R}^{3\times2}} \frac{|g(\tilde{G})|}{|\tilde{G}|^2 + 1} < +\infty \right\},  
\end{align}
and for every measurable $A \subset \omega$ whenever the sequence $\{ f(\nabla y_j) \}$ converges.  As an immediate consequence, we obtain the identities 
\begin{align} \label{eq:YMx}
\langle \nu_x , \id \rangle = \bar{F}, \quad \langle \nu_x, W^{mem} \rangle = W^{mem}(\bar{F}) \quad \text{a.e. } x\in \omega.\end{align}
Now, since $W_{2D}$ is a normal integrand, the fundamental theorem of Young measures gives an inequality, $\liminf_{j \rightarrow \infty} I_{2D}(y_j) \geq  \int_\omega \langle \nu_x, W_{2D} \rangle dx$ (cf. Definition 6.27 and Theorem 8.6, \cite{FL07}).  Thus, 
\begin{align*}
|\omega| W^{mem}(\bar{F}) = \lim_{j \rightarrow \infty} I_{2D}(y_j) \geq 
 \int_\omega \langle \nu_x, W_{2D} \rangle  dx \geq \int_\omega \langle \nu_x, W^{mem} \rangle dx = |\omega| W^{mem}(\bar{F}) 
\end{align*}
where we use the fact that $W_{2D} \geq W^{mem}$.  
It follows $|\omega| W^{mem}(\bar{F}) =  \int_\omega \langle \nu_x, W_{2D}  \rangle dx$.  
Again using the fact that $W_{2D} \geq W^{mem}$ and (\ref{eq:YMx}) we conclude
\begin{align}
\langle \nu_x, W_{2D}  \rangle = W^{mem}(\bar{F})  \quad \text{a.e. } x \in \omega.
\end{align}
By the localizing  properties of $H^1$ gradient Young measures (cf. Theorem 2.3 of \cite{KP91}), we conclude that the fine-scale features which arise from minimizing sequences of $W_{2D}$ are described by the homogenous $H^1$ gradient Young measures which admit the identities,
\begin{align} \label{eq:YMpre}
\langle \nu, \id \rangle = \bar{F}, \;\;\;  \langle \nu, W_{2D} \rangle = W^{mem}(\bar{F}).  
\end{align}
\vspace{\baselineskip}
We present a characterization of this in the following theorem.

\begin{theorem} \label{th:char}
Let $r > 1$, $\bar{F} \in \mathbb{R}^{3\times2}$ and let 
\begin{align} \label{eq:YMWmem}
M_{\bar{F}} := \left\{ \nu \in M^{qc}, \  \langle \nu, \id \rangle = \bar{F}, \ \langle \nu, W_{2D} \rangle = W^{mem}(\bar{F}) \right\} 
\end{align}
be the set of homogenous $H^1$ gradient Young measures that satisfy (\ref{eq:YMpre}).  Then, there exists $\bar{\nu} \in M_{\bar{F}}$.  Further the following is true.

\begin{enumerate}

\item \label{it:m} (The region $\mathcal M$)  Suppose $(\lambda_M(\bar{F}), \delta(\bar{F})) \in \mathcal{M}$.   Set  $\bar{\delta} := \delta(\bar{F})$.  Let the singular value decomposition (cf. (\ref{eq:pd})) of $\bar{F}$ be given by $\bar{F} = \bar{Q} D_{\bar{F}} \bar{R}$ with $\bar{Q} \in SO(3), \bar{R} \in O(2)$.

If $\bar{\nu} \in M_{\bar{F}}$, then
\begin{align}\label{eq:suppnu1}
\supp  \bar{\nu} \subset \mathcal{K}_{\bar{\delta}} :=  \{ \tilde{G} \in \mathbb{R}^{3\times2} \colon
\tilde{G} = \bar{Q} Q D_{\bar{\delta}} R, Q \in SO(3), R \in O(2), \det(R) Q f_3 = \det(\bar{R}) f_3 \}
\end{align}
where $f_3 \in \mathbb{S}^2$ is orthogonal to the plane of the reference configuration of the membrane and
\begin{align}\label{eq:tDTh}
D_{\bar{\delta}} = \bar{\delta}^{1/2} \left( \begin{array}{cc} 
r^{1/4}  & 0\\
0 & r^{-1/4}  \\
0 & 0 
\end{array} \right) .
\end{align}

\item \label{it:w} (The region $\mathcal W$)  Suppose $(\lambda_M(\bar{F}), \delta(\bar{F})) \in \mathcal{W}$. Set $\bar{\lambda}_M = \lambda_M(\bar{F})$.  Further, set $e_M \in \mathbb{S}^2$ and $f_M \in \mathbb{S}^1$ to be the unique pair (up to a change in sign)  of vectors which satisfy $\bar{F} f_M = \bar{\lambda}_M e_M$.  

If $\bar{\nu} \in M_{\bar{F}}$, then
\begin{align*}
\supp \bar{\nu} \subset  \mathcal{K}_{\bar{\lambda}_M} := \{ \tilde{G} \in \mathbb{R}^{3\times2} : ( \lambda_M, \delta)(\tilde{G}) = (\bar{\lambda}_M , \bar{\lambda}_M^{1/2} ) , \;\; \tilde{G} f_M = \bar{\lambda}_M e_M \}.  
\end{align*}

\item \label{it:l} (The region $\mathcal L$)  Suppose $(\lambda_M(\bar{F}), \delta(\bar{F})) \in \mathcal{L}$.

If $\bar{\nu} \in M_{\bar{F}}$, then
\begin{align*}
\supp \bar{\nu} \subset  \mathcal{K}_{0} := \{ \tilde{G} \in \mathbb{R}^{3\times2} : ( \lambda_M, \delta)(\tilde{G}) \in \mathcal{L} \text{ with } \delta/\lambda_M  = r^{-1/6}  \}.  
\end{align*}

\item \label{it:s} (The region $\mathcal S$)  Suppose $(\lambda_M(\bar{F}), \delta(\bar{F})) \in \mathcal{S}$. 
If $\bar{\nu} \in M_{\bar{F}}$, then $\bar{\nu}$ is a Dirac mass.  So $\supp \bar{\nu} = \{ \bar{F} \}$.

\end{enumerate}

\end{theorem}

Theorem \ref{th:char} has striking physical implications.  First consider Part \ref{it:m} corresponding to region $\mathcal{M}$ and consider the particular case when $\bar{Q} = I$.  Consider any $\tilde{G} \in \mbox{supp } \bar{\nu}$ and its characterization in ({\ref{eq:suppnu1}).  Since $\det(R)Q f_3 = \det(\bar{R})f_3$,  it follows $Q D_{\bar{\delta}} R v \cdot f_3 = 0$ for each $v \in {\mathbb R}^2$.  In other words, $Q D_{\bar{\delta}} R$ maps ${\mathbb R}^2$ to ${\mathbb R}^2$.   Thus, all the oscillations are in the plane.   Further, for such matrices $\tilde{G}$,
$$
W_{2D} (\tilde{G}) = W_{0} (\tilde{G}|c) = W^e(\tilde{G}|c,n)
$$
for $c = (0,0,\bar{\delta}^{-1})^T$ and $n\cdot f_3 = 0$.  The first of these identities follows from the fact that 
$(\lambda_M(\tilde{G}), \delta(\tilde{G})) \in \mathcal{S}$ (see Lemma \ref{SupportLemma} below) and Lemma \ref{lem:W2D}, while the second follows from the fact that the largest principal value of $(\tilde{G}|c)$ is $\lambda_M(\tilde{G})$.   Importantly, the director is always in the plane.   In summary, the director oscillates in the plane and oscillations create no out of plane deformation.  The case $\bar{Q} \ne I$ is similar except the plane is oriented by the rotation $\bar{Q}$.  Thus, the fine-scale features in $\mathcal{M}$ is limited to in-plane oscillations of the director.

Now consider Part \ref{it:w}.  First consider the case when $e_M = 
\left(\begin{array}{c} f_M \\0 \end{array}\right)$. Using an argument as before, for any $\tilde{G} \in \mbox{supp } \bar{\nu}$, 
$$
W_{2D} (\tilde{G}) = W_{0} (\tilde{G}|c) = W^e(\tilde{G}|c,n)
$$
for $c \cdot e_M = 0$, $|c| = \bar{\lambda}_M^{-1/2}$ and $n = e_M$.  In other words, the director $n$ is fixed with an in-plane direction $e_M$.  Further, notice $\tilde{G}$ is necessarily of the form
$$
\tilde{G} = Q \left(\begin{array}{cc} \bar{\lambda}_M& 0 \\ 0 & \bar{\lambda}_M^{-1/2} \\ 0&0 \end{array} \right)
$$
in the $e_M - \left( \begin{array}{c} f_M^\perp \\ 0 \end{array} \right)- f_3$ frame for $Q \in SO(3)$ that satisfies $Qe_M = e_M$.  In other words, the membrane is uniformly deformed and the fine features are related to rotations about a fixed axis $e_M$.  In other words, oscillations represent wrinkling and these oscillations are always perpendicular to $e_M$.  The general  case $e_M \ne 
\left(\begin{array}{c} f_M \\0 \end{array}\right)$ is similar except a uniform rotation orients $\left(\begin{array}{c} f_M \\0 \end{array}\right)$ to $e_M$.

Part \ref{it:l} says that region $\mathcal{L}$ involves only the spontaneously deformed states  while Part \ref{it:s} says that there are no fine-scale features in $\mathcal{S}$.
\bigskip

We now turn to the proofs of the theorems.  They rely on the following lemmas.

\begin{lemma}\label{SupportLemma}
Let $\bar{F} \in \mathbb{R}^{3\times2}$ and $\bar{\delta}$ satisfy the hypotheses in Theorem \ref{th:char} Part \ref{it:m}.  Then any $\bar{\nu} \in M_{\bar{F}}$ satisfies
\begin{align}\label{eq:suppYM}
\supp \bar{\nu} \subset  \{ \tilde{G} \in \mathbb{R}^{3\times2} : (\lambda_M,\delta )(\tilde{G}) = (r^{1/4}\bar{\delta}^{1/2},\bar{\delta}) \}.  
\end{align}
\end{lemma}
\begin{lemma}\label{SupportLemma2}
Let $\bar{F} \in \mathbb{R}^{3\times2}$ and $\bar{\lambda}_M$ satisfy the hypotheses of Theorem \ref{th:char} Part \ref{it:w}.  Then any $\bar{\nu} \in M_{\bar{F}}$ satisfies
\begin{align}\label{eq:supportWrinkling}
\supp \bar{\nu} \subset \{ \tilde{G} \in \mathbb{R}^{3\times2} : (\lambda_M, \delta )(\tilde{G}) = (\bar{\lambda}_M,\bar{\lambda}_M^{1/2} )\}.  
\end{align}
\end{lemma}

\noindent {\it Proof of Lemma \ref{SupportLemma}.}  Recall from Section \ref{EffEn} that we may write $W^{mem} = \psi \circ (\lambda_M,\delta)$ and $W_{2D} = \varphi \circ (\lambda_M, \delta)$ where $\psi \;(\varphi) :\mathcal{R} \rightarrow \mathbb{R} \;( \mathbb{R} \cup \{ +\infty\})$ respectively for  $\mathcal{R} = \{(s,t) \in \mathbb{R}^2 : s^2 \geq t, t \geq 0\}$.  Recall also that $\psi$ is a convex, and it is non-decreasing in each argument.  Also, $\psi \leq \varphi$.  Finally, $(\lambda_M, \delta) : \mathbb{R}^{3\times2} \rightarrow \mathcal{R}$ are quasiconvex functions bounded quadratically.  Therefore, for every homogenous  $H^1$ gradient Young measure with $\langle \nu, \id \rangle = \bar{F}$, 
\begin{align}\label{eq:inequalitiesWmem}
W^{mem}(\bar{F}) & = \psi \circ (\lambda_M, \delta) (\langle \nu, \id \rangle) \nonumber \\
& \leq \psi ( \langle \nu , \lambda_M \rangle , \langle \nu, \delta \rangle)  \nonumber \\ 
& \leq \langle \nu, \psi \circ (\lambda_M, \delta) \rangle \nonumber \\
& \leq \langle \nu, \varphi \circ (\lambda_M,\delta) \rangle = \langle \nu, W_{2D} \rangle.  
\end{align}
Here, the first inequality follows from the Jensen's inequality satisfied by homogenous  $H^1$ gradient Young measures  since $(\lambda_M,\delta)$ are quasiconvex with the appropriate growth and $\psi$ is non-decreasing in each argument.  The second inequality follows from the convexity of $\psi$, and the third follows since $\psi \leq \varphi$.  

Now, for any $\bar{\nu} \in M_{\bar{F}}$,  each inequality in  (\ref{eq:inequalitiesWmem}) is an equality.  This restricts the support of $\bar{\nu}$.  To deduce this restriction, suppose that the point $(\lambda_M(\bar{F}), \delta(\bar{F})) \in \mathcal{M}$ corresponds to point C in Figure \ref{fig:a}.

Consider the first inequality.  By quasiconvexity and growth conditions, $\langle \bar{\nu} , \lambda_M \rangle \geq \lambda_M(\bar{F})$ and $\langle \bar{\nu} ,\delta \rangle \geq \delta(\bar{F})$.  In the $\lambda_M-\delta$ space  in Figure \ref{fig:a}, these inequalities imply the point $(\langle \bar{\nu}, \lambda_M \rangle, \langle \bar{\nu}, \delta \rangle)$ cannot be to the left or below point C.  Further, every point to the right and above the point C has higher $\psi$ (cf. Figure \ref{fig:eff}) except the line between and including the points C and D.    Hence,
\begin{align}\label{eq:CD}
(\langle \bar{\nu}, \lambda_M \rangle, \langle \bar{\nu} , \delta \rangle ) \in CD.
\end{align}

 Next, consider the last inequality. Since $\varphi = \psi$ only on $\mathcal{S} \cup \{ (s,t) \in \mathcal{L} : t/s = r^{-1/6}\} =:\mathcal{S}'$  (see Figure \ref{fig:lam}b), we conclude
\begin{align}\label{eq:support} 
\supp \bar{\nu} \subset \Bigl\{ \tilde{G}\in \mathbb{R}^{3\times2} : (\lambda_M,\delta)(\tilde{G}) \in \mathcal{S}' \Bigr\}. 
\end{align}

It remains to consider the middle inequality in (\ref{eq:inequalitiesWmem}).  We do this in Proposition \ref{suppProp} below.  If the middle inequality is an equality, we show in the proposition the support of $\bar{\nu}$ satisfies (\ref{eq:suppYM}).  This completes the proof.  

\endproof

\begin{proposition}\label{suppProp}
Let $\bar{F}$ and $\bar{\delta}$ be as in the Theorem \ref{th:char} Part \ref{it:m}. If $\bar{\nu}$ satisfies (\ref{eq:CD}),(\ref{eq:support}) and
\begin{align} \label{eq:equals}
\psi ( \langle \bar{\nu} , \lambda_M \rangle , \langle \bar{\nu}, \delta \rangle)  =  \langle \bar{\nu}, \psi \circ (\lambda_M, \delta) \rangle,
\end{align}
then the support of $\bar{\nu}$ satisfies (\ref{eq:suppYM}) in Lemma \ref{SupportLemma}. 
\end{proposition}

\noindent {\it Proof.}  Set $\mathcal{A}^+ = \{ \tilde G: (\lambda_M, \delta)(\tilde G) \in
\mathcal{S} \cap \{\delta > \bar{\delta}\} \}$  and
\begin{align*}
\theta^+ = \int_{\mathcal{A}^+} d\bar{\nu}(G) .
\end{align*} 
If $\theta^+ = 1$, then by the polyconvexity of $\delta$, $\langle \bar{\nu}, \delta \rangle > \bar{\delta}$ contradicting (\ref{eq:CD}).  Now consider the case $1 >  \theta^+ > 0 $.  Set
\begin{align*} 
\lambda_M^+ &:= \frac{1}{\theta^+ } \int_{\mathcal{A}^+} \lambda_M(\tilde{G}) d \bar{\nu}(\tilde{G}), 
& \lambda_M^- &:= \frac{1}{1-\theta^+ } \int_{{\mathbb R}^{3\times 2} \setminus \mathcal{A}^+} \lambda_M(\tilde{G}) d \bar{\nu}(\tilde{G}),
\\
\delta^+ &:= \frac{1}{\theta^+ } \int_{\mathcal{A}^+} \delta(\tilde{G}) d \bar{\nu}(\tilde{G}),
& \delta^- &:= \frac{1}{1-\theta^+ } \int_{{\mathbb R}^{3\times 2} \setminus \mathcal{A}^+} \delta(\tilde{G}) d \bar{\nu}(\tilde{G}).
\end{align*}
Clearly, $\delta^+ > \bar{\delta}$  and 
\begin{align}\label{eq:aveLines}
&\theta^+ \lambda_M^+ + (1- \theta^+ ) \lambda_M^- = \langle \bar{\nu}, \lambda_M \rangle, \nonumber \\
&\theta^+ \delta^+ + (1- \theta^+ ) \delta^- = \langle \bar{\nu}, \delta \rangle.
\end{align}
From the equality in (\ref{eq:aveLines}), $\delta^{-} < \bar{\delta}$.  Further, notice from the convexity of $\psi$ that
\begin{align} \label{psiconvex}
\psi (\lambda_M^+, \delta^+) & \le \frac{1}{\theta^+ }  \int_{\mathcal{A}^+}  \psi \left(\lambda_M(\tilde G), \delta(\tilde{G})  \right)d \bar{\nu}(\tilde{G}), \\
\psi (\lambda_M^-, \delta^-) & \le \frac{1}{1-\theta^+ }  \int_{{\mathbb R}^{3\times2} \setminus \mathcal{A}^+}  \psi \left(\lambda_M(\tilde G), \delta(\tilde{G})  \right)d \bar{\nu}(\tilde{G}) \ .
\end{align}

Now, in the $\lambda_M-\delta$ space shown in Figure \ref{fig:a}, the definitions above imply that the point 
$(\lambda_M^+, \delta^+)$ is a point above the line CD while $(\lambda_M^-, \delta^-)$ is below the line CD such that the line joining these points intersect CD.  It is easy to verify by explicitly computing the derivative along such lines (or by inspecting Figure \ref{fig:eff}), that $\psi$ is strictly convex in such segments.    Therefore, 
\begin{align}
\psi( \langle \bar{\nu} , \lambda_M \rangle, \langle \bar{\nu}, \delta \rangle )  
&= \psi \left(\theta^+ \lambda_M^+ + (1- \theta^+ ) \lambda_M^-,\theta^+ \delta^+ + (1- \theta^+ ) \delta^-\right) \nonumber\\
&< \theta^+ \psi(\lambda_M^+, \delta^+)  + (1- \theta^+ ) \psi(\lambda_M^- , \delta^-)\\  
&\leq  \langle \bar{\nu}, \psi \circ (\lambda_M, \delta) \rangle . \nonumber
\end{align}
The last inequality follows from (\ref{psiconvex}).  However, this contradicts the assumption (\ref{eq:equals}).

Therefore, $\theta^+ =0$, and 
\begin{align}\label{eq:support1}
\supp \bar{\nu} \subset \{ (\lambda_M, \delta)(\tilde{G}) \in \mathcal{S}' : \delta(\tilde{G}) \leq \bar{\delta}\},
\end{align}  
which is the compliment of $\mathcal{A}^{+}$ in the set given in (\ref{eq:support}).

Finally, given (\ref{eq:support1}) and since $\bar{\delta} = \langle \bar{\nu}, \delta \rangle$ (see \ref{eq:CD}), it follows that $\supp \bar{\nu} \subset \{ (\lambda_M,\delta)(\tilde{G}) \in \mathcal{S}' : \delta(\tilde{G}) = \bar{\delta} \}$.  But this is just a single point in the $\lambda_M - \delta$ space, and it's given by (\ref{eq:suppYM}).  Thus, we conclude the proposition.  

\endproof

The proof of Lemma \ref{SupportLemma2} is very similar and omitted.
\bigskip

\noindent {\it Proof of Theorem \ref{th:char}.}
Existence of a $\bar{\nu} \in M_{\bar{F}}$ follows from the construction in Section \ref{sec:lb}.

\noindent{\it Part \ref{it:m}.}
For any $\bar{F}$ with $(\lambda_M(\bar{F}), \delta(\bar{F})) \in\mathcal{M}$ and for any $\bar{\nu} \in M_{\bar{F}}$, the support of $\bar{\nu}$ satisfies (\ref{eq:suppYM}) by Lemma \ref{SupportLemma}.   Note that $(\lambda_M,\delta )(\tilde{G}) = (r^{1/4}\bar{\delta}^{1/2}, \bar{\delta})$ is equivalent to stating that the principal values of $\tilde{G}$ are 
$r^{1/4}\bar{\delta}^{1/2}$ and $r^{-1/4}\bar{\delta}^{1/2}$.  Therefore, by the singular value decomposition theorem (\ref{eq:pd}), it follows that
\begin{align}\label{eq:supportSet}
\supp \bar{\nu} \subset \{ \tilde{G} \in \mathbb{R}^{3\times2} : \tilde{G} = Q D_{\bar{\delta}} R, Q\in SO(3), R\in O(2) \} =: \mathcal{K}_{supp}
\end{align}
for $D_{\bar{\delta}}$ is given in (\ref{eq:tDTh}).

Now, for any $D \in {\mathbb R}^{3 \times 2}, Q \in SO(3), R \in O(2)$, it is an easy calculation to find that 
$\adj (QDR) = \det(R) Q \adj D$.   Further for $D$ of the form (\ref{eq:d}), $\adj D = \lambda_1\lambda_2 f_3$. 
Further, the adjugate is a minor and therefore $\langle \bar{\nu} , \adj \rangle = \adj( \langle \bar{\nu}, \id \rangle ) = \adj \bar{F}$.   Recalling the support  (\ref{eq:supportSet}) of $\bar{\nu}$, we conclude 
\begin{align*}
\det(\bar{R}) \bar{Q}f_3 
=\frac{1}{\bar{\delta}}\langle \bar{\nu}, \adj \rangle 
= \frac{1}{\bar{\delta}} \int_{\mathbb{R}^{3\times2}} \adj \tilde{G} d\bar{\nu} (\tilde{G}) 
=  \int_{\mathcal{K}_{supp}} \det (R(\tilde{G})) Q (\tilde{G}) f_3 d \bar{\nu}(\tilde{G}).  
\end{align*}
Note that $\det(\bar{R}) \bar{Q}f_3 \in {\mathbb S}^2$, and $\det (R(\tilde{G})) Q (\tilde{G}) f_3 \in {\mathbb S}^2$ for each $\tilde{G} \in {\mathcal K}_{supp}$.    In other words, the equation above states that an average of a distribution on ${\mathbb S}^2$ yields an element of ${\mathbb S}^2$.  However, since each element of ${\mathbb S}^2$ is an extreme point, it means that the distribution is concentrated at a single point  on $\mathbb{S}^2$.  That is, if we let $Q_{0}(\tilde{G}) = \bar{Q}^T Q(\tilde{G})$, then $\det(\bar{R}) f_3 = \det(R(\tilde{G})) Q_0(\tilde{G}) f_3$. The result follows.

\noindent {\it  Part \ref{it:w}.}
For any $\bar{F}$  with $(\lambda_M(\bar{F}), \delta(\bar{F})) \in\mathcal{W}$ and for any $\bar{\nu} \in M_{\bar{F}}$, it follows from the definition of $e_M, f_M$ that 
\begin{align} \label{eq:average}
\int_{\mathbb{R}^{3\times2}} \tilde{G} f_M  d\bar{\nu}(\tilde{G})  = \bar{\lambda}_M e_M.
\end{align}
So,
\begin{align*}
\int_{\mathbb{R}^{3\times2}} |\tilde{G} f_M|  d\bar{\nu}(\tilde{G})  \ge 
\left|\int_{\mathbb{R}^{3\times2}} \tilde{G} f_M  d\bar{\nu}(\tilde{G})  \right| 
= |\bar{F} f_M| = \bar{\lambda}_M.
\end{align*}
However, from Lemma \ref{SupportLemma2}, we see that 
$\max_{e \in {\mathbb S}^1} |\tilde{G} e| = \lambda_M ( \tilde{G}) = \bar{\lambda}_M$
for each $\tilde{G} \in \mbox{supp } \bar{\nu}$.     Therefore, 
$|\tilde{G} f_M|  \le \bar{\lambda}_M$ for each $\tilde{G} \in \mbox{supp } \bar{\nu}$. 
We conclude that $|\tilde{G} f_M| = \bar{\lambda}_M$ for each $\tilde{G} \in \mbox{supp } \bar{\nu}$.
Setting $\tilde{G} f_M = \bar{\lambda}_M e(\tilde{G})$ for $e(\tilde{G}) \in {\mathbb S}^2$ and
substituting in (\ref{eq:average}), we conclude that $e (\tilde{G})= e_M$ for 
for each $\tilde{G} \in \mbox{supp } \bar{\nu}$.  The result follows.  

\noindent {\it  Part \ref{it:l}.}
For any $\bar{F}$  with $(\lambda_M(\bar{F}), \delta(\bar{F})) \in\mathcal{L}$, the result follows from the fact that $W_{2D}$ is non-negative and $W_{2D}(\tilde G) = 0$ if and only if $\tilde G \in \mathcal{K}_0$. 

\noindent {\it  Part \ref{it:s}.}
Finally, let $\bar{F}  \in \mathbb{R}^{3\times2}$ such that $(\lambda_M(\bar{F}), \delta(\bar{F})) \in \mathcal{S}$, $\bar{\nu} \in M_{\bar{F}}$. Recall $W_{2D} = \varphi \circ (\lambda_M,\delta)$ and $\varphi$ is strictly convex in $\mathcal{S}$.  Thus, 
\begin{align*}
\supp \bar{\nu} \subset \{ \tilde{G} \in \mathbb{R}^{3 \times 2} \colon (\lambda_M, \delta)(\tilde{G}) = (\lambda_M, \delta)(\bar{F}) \}.  
\end{align*}
This is actually equivalent to the set (\ref{eq:DiracMass}) given in Proposition \ref{DiracMass} below since $\lambda_m = \delta/\lambda_M$.  The result follows from the proposition.
\endproof

\begin{proposition}\label{DiracMass}
Let $\bar{F} \in \mathbb{R}^{3 \times 2}$ such that the singular values satisfy the strict inequality
\begin{align*}
\lambda_M(\bar{F}) > \lambda_m(\bar{F})\geq 0.  
\end{align*}
Suppose $\nu$ is a probability measure on the space of $\mathbb{R}^{3\times2}$ matrices such that $\langle \nu, \id \rangle = \bar{F}$ and 
\begin{align}\label{eq:DiracMass}
\supp \nu \subset \{ \tilde{G} \in \mathbb{R}^{3\times2} \colon (\lambda_M , \lambda_m)(\tilde{G}) = (\lambda_M, \lambda_m)(\bar{F})\}.
\end{align}
Then $\nu$ (up to a set of measure zero) is a Dirac mass at $\bar{F}$.
\end{proposition}

\noindent {\it Proof. }To begin, set $(\bar{\lambda}_M , \bar{\lambda}_m) = (\lambda_M(\bar{F}), \lambda_m(\bar{F}))$.  We let $\{ \bar{e}_1, \bar{e}_2\} \subset \mathbb{R}^3$ and $\{ \bar{f}_1, \bar{f}_2 \} \subset \mathbb{R}^2$ be sets of orthonormal vectors such that 
\begin{align}\label{eq:changeVar}
\bar{F} = \bar{\lambda}_M \bar{e}_1 \otimes \bar{f}_1 + \bar{\lambda}_m \bar{e}_2 \otimes \bar{f}_2.
\end{align}
Let $\varphi_{\bar{f}_1}(\tilde{G}) : = |\tilde{G} \bar{f}_1|^2$.  This is a convex function.  Therefore, by Jensen's inequality and given $\langle \nu, \id \rangle = \bar{F}$ with $\bar{F}$ satisfying (\ref{eq:changeVar}), 
\begin{align}\label{eq:Jensen}
\langle \nu , \varphi_{\bar{f}_1} \rangle \geq \varphi_{\bar{f}_1}( \bar{F}) = \bar{\lambda}_M^2.
\end{align}
Conversley, applying a similar change of variables (\ref{eq:changeVar}) to the $\tilde{G} \in \supp \nu$, we see 
\begin{align*}
\langle \nu , \varphi_{\bar{f}_1} \rangle &= \int\left( | (\bar{\lambda}_M e_1 \otimes f_1 + \bar{\lambda}_m e_2 \otimes f_2) \bar{f}_1 |^2\right) (\tilde{G}) d \nu (\tilde{G}) \\
&= \int \left(\bar{\lambda}_M^2 \cos(\theta(\tilde{G}))^2 + \bar{\lambda}_m^2 \sin(\theta(\tilde{G}))^2 \right)d \nu (\tilde{G}) \\
&\begin{cases}
= \bar{\lambda}_M^2 &\text{ if } \;\;  \nu( \{ \tilde{G} \in \mathbb{R}^{3\times 2} \colon \sin(\theta(\tilde{G})) \neq 0 \} ) = 0 \\
< \bar{\lambda}_M^2 &\text{ otherwise},
\end{cases}
\end{align*} 
since by assumption $\bar{\lambda}_M > \bar{\lambda}_m$.   Here, $\cos \theta$ denotes the direction cosine between $f_1$ and $\bar{f}_1$.  Combining this observation with (\ref{eq:Jensen}), we deduce (up to a set of measure zero), $\sin(\theta(\tilde{G})) = 0$.   This implies (up to a change in sign) $f_1 = \bar{f}_1$ in measure.  Since $f_1$ and $f_2$ are orthogonal, it follows that (up to a change in sign) $f_2 = \bar{f}_2$ in measure.  

We repeat this argument substituting $\varphi_{\bar{f}_1}$ with the convex function $\varphi_{\bar{e}_1}(\tilde{G}) = |\tilde{G}^T\bar{e}_1|^2$.  It follows that (up to a change in sign) $e_1 = \bar{e}_1$ and $e_2 = \bar{e}_2$ in measure.  The fact that $\langle \nu, \id \rangle = \bar{F}$ ensures the eigenvectors are fixed and not oscillating in sign with some non-zero measure.  The conclusion follows.  
\endproof

\section{State of stress and connection to tension field theory}\label{stress1}

In this section, we seek to understand the state of stress in the membrane.

Formally, consider an incompressible energy density $W_{3D}$ of the form in (\ref{eq:W3D}) and assume $W_0$ is $C^1$ differentiable.  The Piola-Kirchhoff and the Cauchy stress are defined as
\begin{align} \label{stress}
P (F) = \nabla_F W_0 (F) - p (\adj F)^T  \quad \sigma(F) = (\nabla_F W_0 (F)) F^T - p I
\end{align}
where $p$ is the indeterminate pressure (Lagrange multiplier to enforce incompressibility) and $I$ is identity.  We find $p$ by requiring the tractions to be zero on faces of the membrane.  Alternately, recall that we obtain the membrane energy density $W_{2D}$ by writing $F = (\tilde{F}|c)$ and minimizing with $c$ (when $\tilde{F}$ is full rank).   The minimizer $c_{\tilde F}$ satisfies
\begin{align}
\nabla_c W_0 (\tilde{F}|c_{\tilde F}) - p( \adj \tilde{F}) = 0, \quad c_{\tilde F} \cdot \adj \tilde{F} = \det F = 1  \implies p =\nabla_c W_0 (\tilde{F}|c_{\tilde F}) \cdot c_{\tilde F}. 
\end{align}
Above, $\nabla_c$ denotes derivative with respect to the third column of the deformation gradient.
Substituting this back in (\ref{stress}) and writing $\nabla_F W_0 = (\nabla_{\tilde F} W_0 | \nabla_c W_0)$  we obtain a characterization of the state of stress in the membrane.
\begin{align}\label{eq:Cauchy0}
P_{2D} (\tilde F) &:= P (\tilde F|c_{\tilde F}) = (\nabla_{\tilde{F}} W_0| 0) = (\nabla_{\tilde{F}} W_{2D}|0), \nonumber \\
 \sigma_{2D}(\tilde F) &:= \sigma (\tilde F|c_{\tilde F}) = (\nabla_{\tilde{F}} W_{2D}) {\tilde F}^T.
\end{align}
Notice that these depend only on $W_{2D}$.

However, the effective energy of the membrane is not $W_{2D}$ but its relaxation.  In other words, energy minimization with the integral of $W_{2D}$ can lead to fine-scale oscillations, and thus the stress may also oscillate on a fine scale.  Therefore, we need to understand the overall of effective stress.  Ball et al. \cite{BKK00} have shown that if $f: \mathbb{R}^{n\times m} \rightarrow \mathbb{R}$ is differentiable and satisfies certain growth conditions, then $f^{qc}$ is a $C^1$ function.  Moreover, $\nabla f^{qc}$ can be written in terms of $\nabla f$ and a  homogeneous $W^{1,p}$ gradient Young measure $\nu$ generated by minimizing sequences of $\int_{\Omega} f(\nabla y)dx$, i.e.
\begin{align}\label{eq:gradfqc}
\nabla f^{qc} = \int \nabla f d\nu . 
\end{align}
Unfortunately $W_{2D}$ is an extended function (equal to $+\infty$ when $\rank \tilde F < 2$), and the analogous result is unknown.  However,  our resulting effective energy $W_{2D}^{qc} \equiv W^{mem}$ is finite  everywhere and is differentiable except on a boundary.  So have the following characterization of the stress.

\begin{theorem}\label{StressTheorem}
Let $r > 1$, let $\mathcal{D}\subset \mathbb{R}^2$ be the open set $\mathcal{D} := \{ (s,t) \in \mathbb{R}^2_+ \colon 0<t < s^2 \}$, and let $\bar{F} \in \mathbb{R}^{3\times2}$ such that $(\lambda_M(\bar{F}) , \delta(\bar{F})) \in \mathcal{D}$.  If $\nu_{\bar{F}}$ is a homogenous $H^1$ gradient Young measure generated by minimizing sequences for the energy $I_{2D}$ in the space $\mathcal{A}_{\bar{F}}$ (see (\ref{eq:I2D})) with support in $\mathcal{D}$, then
\begin{align}
& \nabla_{\bar{F}} W^{mem} = \int \nabla_{\tilde{G}} W_{2D} d \nu_{\bar{F}}(\tilde{G}), \label{eq:W2DIdent} \\
& (\nabla_{\bar{F}} W^{mem} )\bar{F}^T = \int (\nabla_{\tilde{G}} W_{2D}) \tilde{G}^T d \nu_{\bar{F}}(\tilde{G}).  \label{CauchyTheorem}
\end{align}
Further, the Cauchy stress $\sigma^{mem}(\bar{F}) := (\nabla_{\bar{F}} W^{mem} )\bar{F}^T $ has the following explicit characterization.
\begin{align}\label{eq:memStress}
\sigma^{mem}
&=\mu r^{1/3} \begin{dcases} 
0 &\text{ if } (\lambda_M, \delta) \in \mathcal{L} \cap \mathcal{D}, \\
\left(\frac{\delta}{r^{1/2}} - \frac{1}{\delta^2}\right) \text{Id}_2 &\text{ if } (\lambda_M, \delta) \in \mathcal{M} \cap \mathcal{D}, \\
\left(\frac{\lambda_M^2}{r} - \frac{1}{\lambda_M} \right) e_1 \otimes e_1 &\text{ if } (\lambda_M,\delta) \in \mathcal{W}  \cap \mathcal{D}, \\
\left(\frac{\lambda_M^2}{r} - \frac{1}{\delta^2}\right) e_1 \otimes e_1 + \left( \frac{\delta^2}{\lambda_M^2} - \frac{1}{\delta^2}\right) e_2 \otimes e_2 &\text{ if } (\lambda_M, \delta) \in \mathcal{S} \cap \mathcal{D}.
\end{dcases}
\end{align}
\end{theorem}

\bigskip

\noindent Before we prove the theorem, we make a few comments on the physical implications.  First, the membrane is always in a state of plane stress in the tangent plane.    Second, the principal stresses (the eigenvalues) are always non-negative.  Therefore, the membrane can not sustain compressive stress.  Further, the stress is zero in region $\mathcal{L}$, uniaxial tension in $\mathcal{W}$, equi-biaxial tension in $\mathcal{M}$ and biaxial tension in $\mathcal{S}$.  
The different regimes are shown in Figure  \ref{fig:stress}.

To understand this further, consider the special case $r=1$ when this theory reduces to that of neo-Hookean elastic membrane.  The region ${\mathcal M}$ now disappears and we are left with regions ${\mathcal L, \mathcal{W}}$ and ${\mathcal S}$ with zero, uniaxial tension and biaxial tension respectively as in the traditional tension field theory \cite{M70,Pip86, SP89}.

Nematic elastomers membranes with $r>1$ are characterized by an additional region ${\mathcal M}$ where the state of stress is equi-biaxial tension.  This is true even though, the principal stretches $(\lambda_M, \delta/\lambda_M)$ can be unequal.  In other words, one can have shear strain but no shear stress.   This is a potentially useful attribute of liquid crystal elastomers in membrane applications. 

\begin{figure}
\centering%
\includegraphics[width=3in]{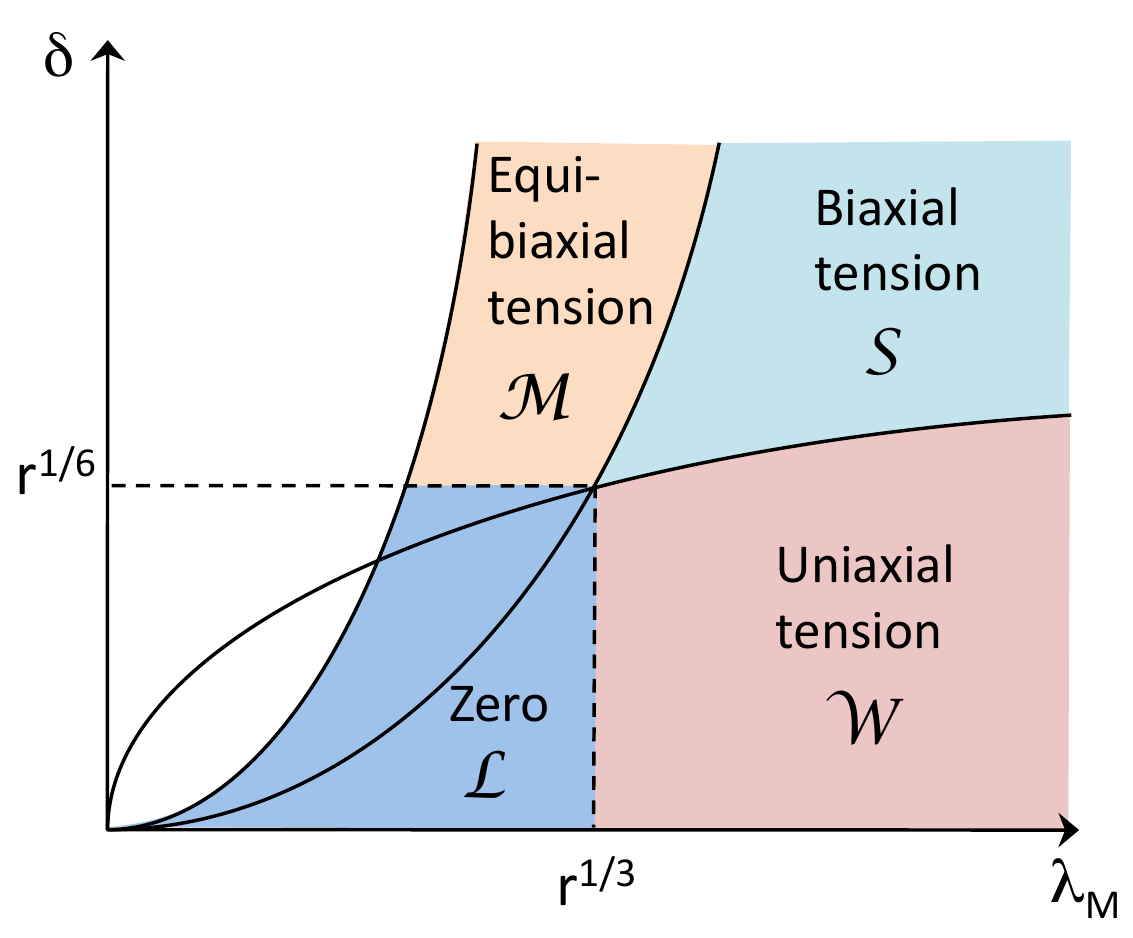}
\caption{The effective stress of nematic elastomer membranes.}\label{fig:stress}
\end{figure}

We turn now to the proof of Theorem \ref{StressTheorem}.
\bigskip

\noindent {\it Proof of Theorem \ref{StressTheorem}.}  Recall from Section \ref{EffEn} that we may write $W_{2D} = \varphi \circ (\lambda_M, \delta)$ and $W^{mem} = \psi \circ (\lambda_M, \delta)$.  
Now,  any $\tilde{G} \in \mathbb{R}^{3\times2}$ has the representation
\begin{align} \label{eq:rep}
\tilde{G} = \lambda_M e_1 \otimes f_1 + \lambda_m e_2 \otimes f_2,
\end{align}
where $\{ e_1, e_2\} \subset \mathbb{R}^3$ and $ \{f_1,f_2\} \subset \mathbb{R}^2$ are orthonormal and $\lambda_M \geq \lambda_m \geq 0$ are the singular values of $\tilde{G}$.  These singular values (and therefore $\delta = \lambda_m\lambda_M$) are continuously differentiable with respect to $\tilde{G}$ as long as they are distinct, i.e. $\lambda_M > \lambda_m$ with
\begin{align*}
\nabla_{\tilde{G}} \lambda_M = e_1 \otimes f_1, \;\;\; \nabla_{\tilde{G}} \lambda_m = e_2 \otimes f_2,
\end{align*}
(cf.  Corollary 3.5 and Theorem 5.1, \cite{QW96}).  We can use this fact and the representation for $\varphi, \psi$ in  Theorem \ref{1407231520} and Lemma \ref{lem:W2D} to conclude that $W_{2D}$ and $W^{mem}$ are continuously differentiable on $\{\tilde{G}: (\lambda_M(\tilde{G}), \delta(\tilde{G}) ) \in \mathcal{D} \} $.  

The rest of the proof is by computation and verification.
\bigskip

\noindent {\it Case 1: $(\lambda_M(\bar{F}), \delta(\bar{F})) \in \mathcal{M} \cap \mathcal{D}$.}  
Set $\bar{\delta} = \delta(\bar{F})$.   According to Theorem \ref{th:char} Part \ref{it:m},  $\supp \nu_{\bar{F}} \subset \mathcal{K}_{\bar{\delta}}$.    We can now apply the representation (\ref{eq:rep}) to $\bar{F}$ and $\tilde{G} \in \supp \bar{\nu}$ to write the identity 
$\bar{F} = \langle \nu_{\bar{F}}, \id \rangle $ as
\begin{align}\label{eq:nuIdent1}
\bar{\lambda}_M \bar{e}_1 \otimes \bar{f}_1 + \left(\frac{\bar{\delta}}{\bar{\lambda}_M} \right) \bar{e}_2 \otimes \bar{f}_2 = \bar{\delta}^{1/2} \int_{\mathcal{K}_{\bar{\delta}}} \left(r^{1/4} e_1 \otimes f_1 + r^{-1/4} e_2 \otimes f_2 \right)(\tilde{G}) d \nu_{\bar{F}}(\tilde{G}).  
\end{align}
Another implication of Theorem \ref{th:char} is any $\tilde{G} \in \supp \nu_{\bar{F}}$ can be written as $\bar{Q} Q D_{\bar{\delta}} R$ where $\bar{Q} \in SO(3)$ arises from the identity $\bar{F} = \bar{Q} D_{\bar{F}} \bar{R}$, for some $Q \in SO(3)$ and $R \in O(2)$ such that $\det(R) Qf_3 = \det(\bar{R})f_3$.  Here, $f_3 \in \mathbb{S}^2$ is orthogonal to the reference configuration of the membrane.  Without loss of generality, we assume $f_3 = (0,0,1)^T$.  Now for each $\alpha = 1, 2$ there is a corresponding $c_\alpha>0$ such that $e_\alpha \cdot (\bar{Q}f_3) = c_\alpha (\bar{Q} Q D_{\bar{\delta}} R f_\alpha) \cdot (\bar{Q}f_3) = c_\alpha 
(D_{\bar{\delta}} R f_\alpha) \cdot (Q^T f_3) = c_{\alpha} (\det R / \det \bar{R}) (D_{\bar{\delta}} R f_{\alpha}) \cdot f_3 = 0$.  In other words, the vectors $e_1$ and $e_2$ span the plane perpendicular to $\bar{Q}f_3$ for each $\tilde{G} \in \supp \nu_{\bar{F}}$. Moreover $\bar{Q}f_3 = \bar{e}_1 \times \bar{e}_2$, and therefore $\bar{e}_1$ and $\bar{e}_2$ also span this plane.  Now, let $R_0$ be a $90$ degree rotation about $f_3$ and $Q_0$ be a $90$ degree rotation about $\bar{Q}f_3$ so that $\bar{R}_0^T \bar{f}_1 = \bar{f}_2,  \bar{R}_0^T \bar{f}_2 = -\bar{f}_1$ and $\bar{Q}_0 \bar{e}_1 = \bar{e}_2, \bar{Q}_0 \bar{e}_2 = -\bar{e}_1$.    Since $e_1-e_2$ span the same plane as $\bar{e}_1-\bar{e}_2$, and $f_1-f_2$ the same plane as $\bar{f}_1-\bar{f}_2$, we have the following relation.  
If $\det(\bar{R}) = 1$,  then 
\begin{align}\label{eq:1Invariant}
&\bar{R}_0^T f_1 = f_2 , \;\; \bar{Q}_0 e_1 = e_2 \;\; \text{ if } \;\; \det(R) = 1, \nonumber \\
&\bar{R}_0^T f_1 = -f_2, \;\; \bar{Q}_0 e_1 = -e_2 \;\; \text{ if } \;\; \det(R) = -1.
\end{align}
If $\det(\bar{R}) = -1$, then
\begin{align}\label{eq:2Invariant}
&\bar{R}_0^T f_1 = -f_2, \;\; \bar{Q}_0 e_1 = -e_2, \;\; \text{ if } \;\; \det(R) = 1, \nonumber \\
&\bar{R}_0^T f_1 = f_2, \;\; \bar{Q}_0 e_1 = e_2, \;\; \text{ if } \;\; \det(R) = -1.  
\end{align}
Thus, pre-multiplying and post-multiplying the identity in (\ref{eq:nuIdent1}) by $\bar{\delta}^{-1/2}\bar{Q}_0$ and $\bar{R}_0$ respectively yields the identity
\begin{align}\label{eq:ForMStress}
\left(\frac{\bar{\lambda}_M}{\bar{\delta}^{1/2}}\right) \bar{e}_2 \otimes \bar{f}_2 + \left( \frac{\bar{\delta}^{1/2}}{\bar{\lambda}_M} \right) \bar{e}_1 \otimes \bar{f}_1 &= \int_{\mathcal{K}_{\bar{\delta}}} \left(r^{1/4} \bar{Q}_0e_1 \otimes \bar{R}_0^T f_1 + r^{-1/4} \bar{Q}_0e_2 \otimes \bar{R}_0^Tf_2 \right)(\tilde{G}) d \nu_{\bar{F}}(\tilde{G}) \nonumber \\
&= \int_{\mathcal{K}_{\bar{\delta}}} \left(r^{1/4} e_2 \otimes f_2 + r^{-1/4} e_1 \otimes f_1 \right)(\tilde{G}) d \nu_{\bar{F}}(\tilde{G})
\end{align}
by (\ref{eq:1Invariant}) and (\ref{eq:2Invariant}).  

It is easy to verify $\frac{\partial \varphi}{\partial \lambda_M}  (r^{1/4} \bar{\delta}^{1/2}, \bar{\delta}) = 0$.  Thus, combining explicit differentiation evaluated in $\mathcal{K}_{\bar{\delta}}$ with the identity (\ref{eq:ForMStress}), we observe
\begin{align}\label{eq:BigCalc}
\int \nabla_{\tilde{G}} W_{2D} d\nu_{\bar{F}} (\tilde{G}) 
&= \int \left(\frac{\partial \varphi}{\partial \lambda_M}  \nabla_{\tilde{G}} \lambda_M 
+ \frac{\partial \varphi}{\partial \delta} \nabla_{\tilde{G}} \delta \right) d \nu_{\bar{F}}(\tilde{G})   \nonumber \\ 
&= \int_{\mathcal{K}_{\bar{\delta}}} \left( \frac{\partial \varphi}{\partial \delta}  \left[ \frac{\delta}{\lambda_M} e_1 \otimes f_1 + \lambda_M e_2 \otimes f_2\right] \right) (\tilde{G}) d \nu_{\bar{F}}(\tilde{G}) \nonumber  \\
&= \mu r^{1/3} \int_{\mathcal{K}_{\bar{\delta}}} \left( \left[ \frac{\delta}{\lambda_M^2} - \frac{1}{\delta^3} \right] \left[ \frac{\delta}{\lambda_M} e_1 \otimes f_1 + \lambda_M e_2 \otimes f_2\right] \right)(\tilde{G}) d \nu_{\bar{F}} (\tilde{G})  \nonumber \\
&= \mu r^{1/3}  \left(\frac{\bar{\delta}^{1/2}}{r^{1/2}} - \frac{1}{\bar{\delta}^{3/2}} \right) \int_{\mathcal{K}_{\bar{\delta}}} \left( r^{-1/4} e_1 \otimes f_1 + r^{1/4} e_2 \otimes f_2 \right)(\tilde{G}) d \nu_{\bar{F}} (\tilde{G}) \nonumber \\
&=  \mu r^{1/3}  \left(\frac{\bar{\delta}^{1/2}}{r^{1/2}} - \frac{1}{\bar{\delta}^{3/2}} \right) \left(\left( \frac{\bar{\delta}^{1/2}}{\bar{\lambda}_M} \right) \bar{e}_1 \otimes \bar{f}_1+ \left(\frac{\bar{\lambda}_M}{\bar{\delta}^{1/2}}\right) \bar{e}_2 \otimes \bar{f}_2 \right)  \nonumber \\
&= \mu r^{1/3} \left\{\left( \frac{\bar{\delta}}{\bar{\lambda}_M r^{1/2} } - \frac{1}{\bar{\lambda}_M \bar{\delta}} \right) \bar{e}_1 \otimes \bar{f}_1+ \left( \frac{\bar{\lambda}_M}{r^{1/2}} - \frac{\bar{\lambda}_M}{\bar{\delta}^2} \right) \bar{e}_2 \otimes \bar{f}_2\right\}.
\end{align}
Finally, it can be verified explicitly that $\nabla_{\bar{F}} W^{mem}$ coincides with (\ref{eq:BigCalc}).    This gives the identity (\ref{eq:W2DIdent}) for region $\mathcal{M}\cap \mathcal{D}$.  

Similarly,
\begin{align}\label{eq:CauchyCalc2}
\int \nabla_{\tilde{G}} W_{2D} \tilde{G}^T d \nu_{\bar{F}}(\tilde{G}) 
&= \int_{\mathcal{K}_{\bar{\delta}}} \left( \frac{\partial \varphi}{\partial \delta}  \nabla_{\tilde{G}} \delta \right) \tilde{G}^T d \nu_{\bar{F}} (\tilde{G}) \nonumber  \\
&= \mu r^{1/3} \int_{\mathcal{K}_{\bar{\delta}}} \left( \frac{\delta^2}{\lambda_M^2} - \frac{1}{\delta^2} \right) \left( e_1 \otimes e_1 + e_2 \otimes e_2\right) (\tilde{G}) d \nu_{\bar{F}} (\tilde{G})  \nonumber  \\
&= \mu r^{1/3} \left( \frac{\bar{\delta}}{r^{1/2}} - \frac{1}{\bar{\delta}^2} \right) \int_{\mathcal{K}_{\bar{\delta}}} (e_1 \otimes e_1 + e_2 \otimes e_2) (\tilde{G}) d \nu_{\bar{F}}(\tilde{G}) \nonumber \\
&= \mu r^{1/3} \left( \frac{\bar{\delta}}{r^{1/2}} - \frac{1}{\bar{\delta}^2} \right) \text{Id}_2.  
\end{align}
The fourth equality uses the fact that the basis $\{ e_1(\tilde{G}) , e_2(\tilde{G})\}$ always spans the same plane.  
Finally, it can be verified explicitly that $\nabla_{\bar{F}} W^{mem} \bar{F}^T$ coincides with (\ref{eq:CauchyCalc2}) in this region.    This gives the identities (\ref{CauchyTheorem}) and (\ref{eq:memStress}) for region $\mathcal{M}\cap \mathcal{D}$.  
\bigskip

\noindent {\it Case 2: $(\lambda_M(\bar{F}), \delta(\bar{F})) \in \mathcal{W} \cap \mathcal{D}$.}  
Set $\bar{\lambda}_M = \lambda_M(\bar{F})$.  Following Theorem \ref{th:char} Part \ref{it:w}, $\supp \nu_{\bar{F}} \subset \mathcal{K}_{\bar{\lambda}_M}$ and so any $\tilde{G} \in \supp \nu_{\bar{F}}$ satisfies $(\lambda_M(\tilde{G}), \delta(\tilde{G})) = (\bar{\lambda}_M, \bar{\lambda}_M^{1/2})$.  In addition, for the vectors $\bar{f}_1 \in \mathbb{S}^1$ and $\bar{e}_1 \in \mathbb{S}^2$ such that $\bar{F} \bar{f}_1 = \bar{e}_1$, $\tilde{G} \in \supp \nu_{\bar{F}}$ also satisfies $\tilde{G} \bar{f}_1 = \bar{e}_1$.  Writing $\tilde{G} \in \supp \nu_{\bar{F}}$ as in (\ref{eq:rep}), we observe using the properties of the set $\mathcal{K}_{\bar{\lambda}_M}$,
\begin{align}\label{eq:WeVec}
\tilde{G} \bar{f}_1 &= (\lambda_M e_1 \otimes f_1 + (\delta / \lambda_M) e_2 \otimes f_2) \bar{f}_1 \nonumber \\
&=( \bar{\lambda}_M e_1 \otimes f_1 + \bar{\lambda}_M^{-1/2} e_2 \otimes f_2) \bar{f}_1 \nonumber \\
&= \bar{\lambda}_M \cos(\theta) e_1 + \bar{\lambda}_M^{-1/2} \sin(\theta) e_2 = \bar{\lambda}_M \bar{e}_1.
\end{align}
Here, $\cos(\theta)$ denotes the direction cosine from $\bar{f}_1$ to $f_1$.  Applying the squared norm to the identities in (\ref{eq:WeVec}) yields $|\tilde{G} \bar{f}_1|^2 = (\bar{\lambda}_M)^2 \cos(\theta)^2 + \bar{\lambda}_M^{-1} \sin(\theta)^2 = \bar{\lambda}_M^2$.  Since $\bar{\lambda}_M^2 > \bar{\lambda}_M^{-1}$ in $\mathcal{W}$, we deduce from this equation that $\cos(\theta) = \pm 1$.  That is, $f_1$ is up to a change in sign equal to $\bar{f}_1$.  Substituting for $f_1$ back into (\ref{eq:WeVec}), we find $e_1 = \pm \bar{e}_1$ when $f_1 = \pm \bar{f}_1$, or alternatively
\begin{align}\label{eq:fixedVecs}
e_1 \otimes f_1 = \bar{e}_1 \otimes \bar{f}_1 \quad \forall \;\; \tilde{G} \in \supp \nu_{\bar{F}}.  
\end{align}

Now, it is easy to verify explicitly $\frac{\partial \varphi}{\partial \delta} ( \bar{\lambda}_M,\bar{\lambda}_M^{1/2}) = 0$.  Thus, combining explicit differentiation evaluated in $\mathcal{K}_{\bar{\lambda}_M}$ with (\ref{eq:fixedVecs}),
\begin{align}\label{eq:Comp1}
\int \nabla_{\tilde{G}} W_{2D} d\nu_{\bar{F}} (\tilde{G}) 
&= \int_{\mathcal{K}_{\bar{\lambda}_M}} \left( \frac{\partial \varphi}{\partial \lambda_M} \nabla_{\tilde{G}} \lambda_M \right) d \nu_{\bar{F}}(\tilde{G}) \nonumber  \\
&= \mu r^{1/3} \int_{\mathcal{K}_{\bar{\lambda}_M}} \left( \left[\frac{\lambda_M}{r} - \frac{\delta^2}{\lambda_M^3} \right] e_1 \otimes f_1\right) (\tilde{G}) d \nu_{\bar{F}}(\tilde{G}) \nonumber \\
&= \mu r^{1/3} \left( \frac{\bar{\lambda}_M}{r} - \frac{1}{\bar{\lambda}_M^2} \right) \int_{\mathcal{K}_{\bar{\lambda}_M} }e_1 \otimes f_1  (\tilde{G}) d \nu_{\bar{F}}(\tilde{G}) \nonumber \\
 &=  \mu r^{1/3} \left(\frac{\bar{\lambda}_M}{r} - \frac{1}{\bar{\lambda}_M^2} \right) \bar{e}_1 \otimes \bar{f}_1.
\end{align}
Finally, it can be verified explicitly that $\nabla_{\bar{F}} W^{mem}$ coincides with (\ref{eq:Comp1}).   Therefore, the identity (\ref{eq:W2DIdent}) is satisfied for $\mathcal{W}$.  

Similarly, 
\begin{align}\label{eq:CauchyCalc1}
\int \nabla_{\tilde{G}} W_{2D} \tilde{G}^T d \nu_{\bar{F}}(\tilde{G}) 
&= \int_{\mathcal{K}_{\bar{\lambda}_M}} \left( \frac{\partial \varphi}{ \partial \lambda_M} \nabla_{\tilde{G}} \lambda_M \right) \tilde{G}^T d \nu_{\bar{F}} (\tilde{G})  \nonumber \\
&=  \mu r^{1/3} \int_{\mathcal{K}_{\bar{\lambda}_M}} \left( \frac{\lambda_M^2}{r^2} - \frac{\delta^2}{\lambda_M^2}\right) e_1 \otimes e_1 (\tilde{G}) d \nu_{\bar{F}}(\tilde{G}) \nonumber \\
&=  \mu r^{1/3} \left( \frac{\bar{\lambda}_M^2}{r} - \frac{1}{\bar{\lambda}_M} \right) \int_{\mathcal{K}_{\bar{\lambda}_M}} e_1 \otimes e_1 (\tilde{G}) d \nu_{\bar{F}}(\tilde{G}) \nonumber \\
&= \mu r^{1/3} \left( \frac{\bar{\lambda}_M^2}{r} - \frac{1}{\bar{\lambda}_M} \right) \bar{e}_1 \otimes \bar{e}_1.
\end{align} 
For the last equality, recall $e_1 = \pm \bar{e}_1$ for $\tilde{G} \in \supp \nu_{\bar{F}}$.  Finally, it is easy to verify explicitly that $\nabla_{\bar{F}} W^{mem} \bar{F}^T$ coincides with (\ref{eq:CauchyCalc1}) in this region.    Thus, we have the identities (\ref{CauchyTheorem}) and (\ref{eq:memStress}) for region $\mathcal{W}\cap \mathcal{D}$.  
\bigskip

\noindent {\it Case 3: $(\lambda_M(\bar{F}), \delta(\bar{F})) \in \mathcal{L} \cap \mathcal{D}$.}  
According to Theorem \ref{th:char} Part \ref{it:l}, $\supp \bar{\nu} \subset \mathcal{K}_0$.  We see that $\nabla W_{2D} = 0$ on $ \mathcal{K}_0$ and similarly $\nabla W^{mem} = 0$ on $ \mathcal{L} \cap \mathcal{D}$.  The identities (\ref{eq:W2DIdent}),  (\ref{CauchyTheorem}) and (\ref{eq:memStress}) for region $\mathcal{L}\cap \mathcal{D}$.  
\bigskip

\noindent {\it Case 4: $(\lambda_M(\bar{F}), \delta(\bar{F})) \in \mathcal{S} \cap \mathcal{D}$.}    
According to Theorem \ref{th:char} Part \ref{it:s}, $\supp \bar{\nu}$ is a Dirac mass.  According to Theorem \ref{1407231520}, $W_{2D}$ and $W^{mem}$ coincide on $\mathcal{S} \cap \mathcal{D}$.   The identities (\ref{eq:W2DIdent}),  (\ref{CauchyTheorem}) and (\ref{eq:memStress}) for region $\mathcal{S}\cap \mathcal{D}$.  
\endproof
}

\subsection*{Acknowledgments}
This work was partially conducted when PC  held a position at the California Institute of Technology. 
We acknowledge support from the US Department of Energy 
National Nuclear Security Administration  (Award Number DE-FC52-08NA28613, all authors),
the US National Science Foundation (Award Number OISE-0967140, PPP and KB)
and the European Research Council under the 
European Union's Seventh Framework Programme (FP7/2007-2013, ERC grant agreement N. 291053, PC). The authors  thank Jan Kristensen for his advice on a draft version of the paper.

\appendix
\section{Appendix A: Proof of Part 2 of Lemma \ref{lem:ub}}


\begin{figure}[h!]
\centering%
\includegraphics[width=11.0cm]{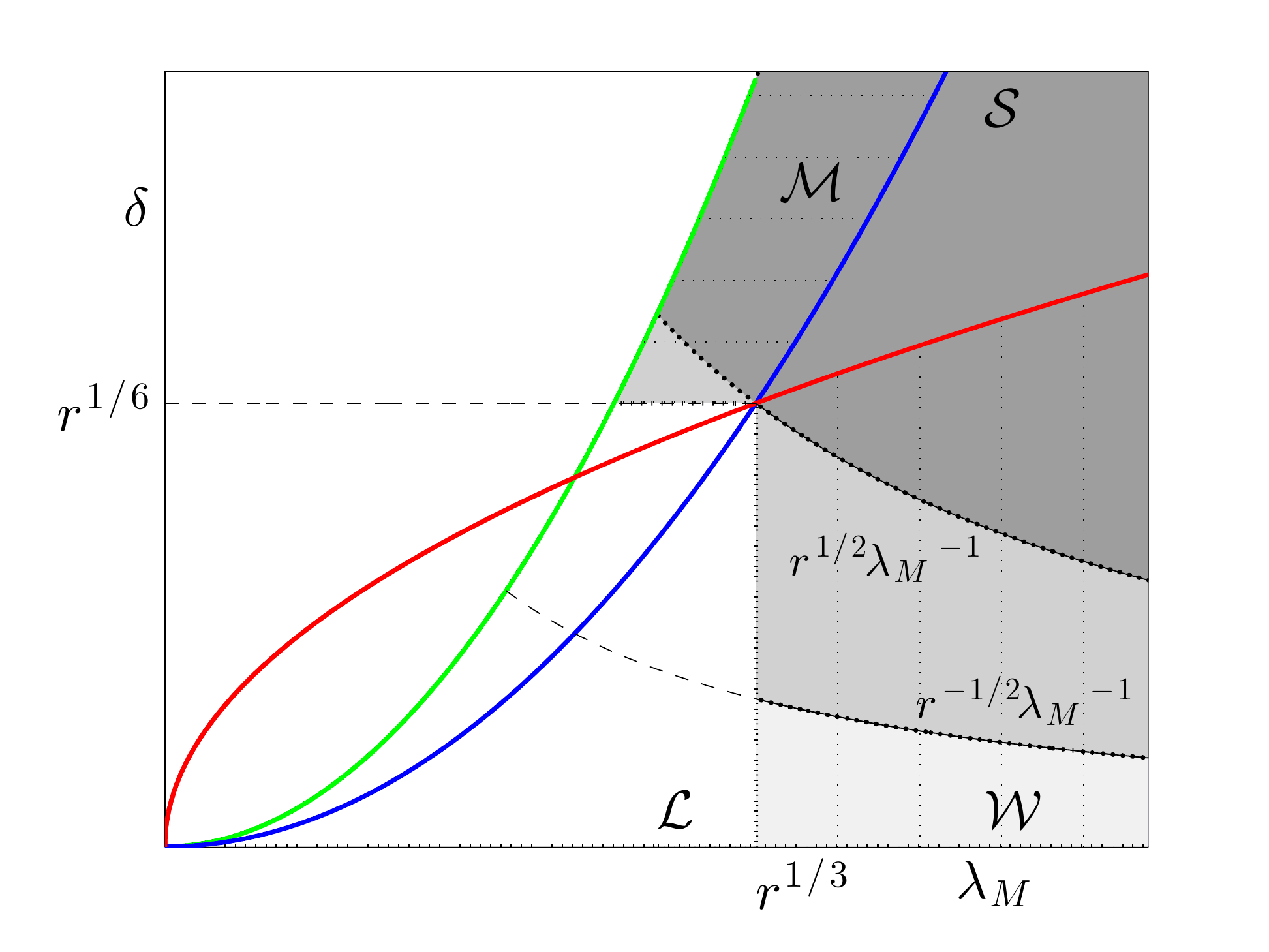}%
\caption{
Regions where the comparison for $W^{mem}$ and the functions $\varphi_1$ (dark gray), $\varphi_2$ (light gray) and $\varphi_3$ (silver/intermediate gray) occurs. The web version of this article contains the above plot figures in color.
}
\label{1309201910}
\end{figure}

\subsection{$W^{mem}\leq \varphi_2$}\label{1309161821}
First, in the region $\mathcal{L}$ of liquid behavior there is nothing to prove. 
Therefore, referring to Fig. \ref{1309201910}, we are left with showing that $W^{mem}\geq  \varphi_2$ in the light gray region of equations for $r^{1/3}< \lambda_M\leq r^{-1/2}\delta^{-1}$.
Recalling that in this region $W^{mem}=\frac{\mu}{2}[r^{1/3} ((\lambda_M(\tilde{F}))^2r^{-1}+2(\lambda_M(\tilde{F}))^{-1} -3]$, it is enough to prove that 
\begin{eqnarray}\label{1303241836}
\frac{\lambda_M^2}{r}+\frac{2}{\lambda_M} \leq
\inf_{\delta} \Bigl\{\Bigl(\lambda_M^2+\frac{\delta^2}{\lambda_M^2}+\frac{1}{r\delta^2}\Bigr) \quad\textrm{for: } \delta\leq r^{-1/2} \lambda_M^{-1}, \,\lambda_M> r^{1/3}\Bigr\}.
\end{eqnarray}
The critical point of $\lambda_M^2+\delta^2\lambda_M^{-2}+r^{-1}\delta^{-2}$ is attained at $\delta^2=\lambda_M r^{-1/2}$. This corresponds to a minimum, yielding the following inequality
\begin{eqnarray}\label{1303241900}
\frac{\lambda_M^2}{r}+\frac{2}{\lambda_M} \leq
\lambda_M^2
+\frac{2}{\lambda_Mr^{1/2}}
\end{eqnarray}
which is indeed true for $\lambda_M>r^{1/3}$. Then, evaluation of  $\lambda_M^2+\delta^2\lambda_M^{-2}+r^{-1}\delta^{-2}$ along the curve $\delta= r^{-1/2} \lambda_M^{-1}$ does not improve the inequality above.

\subsection{$W^{mem}\leq  \varphi_1$}\label{1309161822}
We focus on the interval $\delta\geq r^{1/2}\lambda_M^{-1}$ corresponding (if again we ignore the region $\mathcal{L}$) to the dark gray area in Fig. \ref{1309201910}. This set has a non-empty intersection with both the simple-laminate regions $\mathcal{M}$, $\mathcal{W}$ and the regime of solid behavior $\mathcal{S}$. First of all, 
notice that if $(\lambda_M,\delta)\in \mathcal{S}$ then $W^{mem}\equiv W_{2D}$ and there is nothing to prove.

Let us assume $r^{-1/2}\lambda_M^2<\delta\leq \lambda_M^2$, $\delta\geq r^{1/2}\lambda_M^{-1}$. This corresponds to a subset of $\mathcal{M}$ for which we have 
$W^{mem}=\frac{\mu}{2}[r^{1/3}( 2\delta(\tilde{F})r^{-1/2}+(\delta(\tilde{F}))^{-2} )-3]$. We are left with the inequality
\begin{eqnarray}
\frac{2\delta}{r^{1/2}}+\frac{1}{\delta^2}\leq    
\frac{\lambda_M^2}{r}+\frac{\delta^2}{\lambda_M^2}+\frac{1}{\delta^2}
  \quad\textrm{for }  
r^{-1/2}\lambda_M^2<\delta\leq \lambda_M^2, \delta\geq r^{1/2}\lambda_M^{-1},\nonumber
\end{eqnarray}
which is trivially true.

Then, let us assume $r^{1/2}\lambda_M^{-1}\leq \delta<\lambda_M^{1/2}$. This is a subset of the region $\mathcal{W}$ for which we have 
$W^{mem}=\frac{\mu}{2}\{r^{1/3}[((\lambda_M(\tilde{F}))^2 r^{-1}+2(\lambda_M(\tilde{F}))^{-1}]-3\}$. The inequality
\begin{eqnarray}
\frac{\lambda_M^2}{r}+\frac{2}{\lambda_M}
\leq    
\frac{\lambda_M^2}{r}+\frac{\delta^2}{\lambda_M^2}+\frac{1}{\delta^2}
  \quad\textrm{for }  
 r^{1/2}\lambda_M^{-1}\leq \delta<\lambda_M^{1/2},\nonumber
\end{eqnarray}
follows trivially.

\subsection{$W^{mem}\leq  \varphi_3$}
We now focus on the interval $\delta^{-1}r^{-1/2}<  \lambda_M<\delta^{-1}r^{1/2}$ corresponding to the silver/intermediate gray area in Fig. \ref{1309201910}.
Notice that if we remove the region $\mathcal{L}$ (for which there is nothing to prove), we are left with two disjoint subsets.

We begin with considering $\lambda_M> r^{1/3}$. Since in this region $W^{mem}=\frac{\mu}{2}[r^{1/3} ((\lambda_M(\tilde{F}))^2r^{-1}+2(\lambda_M(\tilde{F}))^{-1} -3]$, it is enough to show that
\begin{eqnarray}\label{}
\frac{\lambda_M^2}{r} +\frac{2}{\lambda_M} 
\leq
\inf\Bigl\{\Bigl(\frac{\delta^2}{\lambda_M^2}+\frac{2\lambda_M}{r^{1/2}\delta} \Bigr)\quad \textrm{for }\frac{1}{\lambda_M}r^{-1/2}<  \delta<\frac{1}{\lambda_M}r^{1/2}, \lambda_M> r^{1/3}\Bigr\},\nonumber
\end{eqnarray}
which is equivalent to
\begin{eqnarray}\label{1309151300}
\frac{\lambda_M^2}{r} +\frac{2}{\lambda_M} 
\leq
\inf\Bigl\{\Bigl(\lambda_m^2+2\frac{1}{r^{1/2}\lambda_m} \Bigr)\quad \textrm{for }\frac{1}{\lambda_M^2}r^{-1/2}<  \lambda_m<\frac{1}{\lambda_M^2}r^{1/2}, \lambda_M> r^{1/3}\Bigr\}.
\end{eqnarray}
The critical point of $\lambda_m^2+2r^{-1/2}\lambda_m^{-1}$ is attained at $(\lambda_m,\lambda_M)=(r^{-1/6},r^{1/3})\in L$. 
Then, we have to evaluate $\lambda_m^2+2r^{-1/2}\lambda_m^{-1}$ on the curves of equations $\lambda_m=\lambda_M^{-2}r^{-1/2}$ and $ \lambda_m=\lambda_M^{-2}r^{1/2}$. 
Notice that, for 
$\lambda_m=r^{1/2} \lambda_M^{-2}$ we have that $\varphi_3=\varphi_1$
while for $\lambda_m=r^{-1/2} \lambda_M^{-2}$ we have $\varphi_3=\varphi_2$ and from discussion of these cases in Paragraph \ref{1309161822} and \ref{1309161821}, respectively, it therefore follows that $(\ref{1309151300})$ is true.

To conclude, we have to prove that the inequality $W^{mem}\leq \varphi_3$ holds in the remaining subregion defined by $\lambda_M^2\geq \delta$, $\delta> r^{1/6}$ and $\delta< r^{1/2}\lambda_M^{-1}$. This subset is contained in the region $\mathcal{M}$ in which case we have $W^{mem}=\frac{\mu}{2}[r^{1/3}( 2\delta(\tilde{F})r^{-1/2}+(\delta(\tilde{F}))^{-2} )-3]$. Therefore, we are left with proving the inequality
\begin{eqnarray}\label{1407211622}
\frac{2\delta}{r^{1/2}}+\frac{1}{\delta^2}\leq \inf\Bigl\{\Bigl(\frac{\delta^2}{\lambda_M^2}+\frac{2\lambda_M}{r^{1/2}\delta}\Bigr) \quad\textrm{for } \delta^{1/2}\leq\lambda_M,\delta> r^{1/6},\delta<\frac{r^{1/2}}{\lambda_M}\Bigr\},\nonumber
\end{eqnarray}
equivalent to
\begin{eqnarray}\label{1407211624}
\frac{2\delta}{r^{1/2}}+\frac{1}{\delta^2}\leq \inf\Bigl\{\Bigl(\lambda_m^2+\frac{2}{r^{1/2}\lambda_m}\Bigr) \quad\textrm{for }\lambda_m\leq \delta^{1/2},\delta> r^{1/6},\lambda_m> r^{-1/2}\delta^2\Bigr\}.
\end{eqnarray}
In order to prove the inequality above it is enough to evaluate the function $\lambda_m^2+2r^{-1/2}\lambda_m^{-1}$ on the boundary of the region defined on the right hand side of $(\ref{1407211624})$.
This yields the following two relations
\begin{eqnarray}\label{1407211722}
\frac{2\delta}{r^{1/2}}+\frac{1}{\delta^2}\leq 
\delta+\frac{2}{r^{1/2}\delta^{1/2}} \quad\textrm{for } \delta\in(r^{1/6},r^{1/3})
\end{eqnarray}
\begin{eqnarray}\label{1407211724}
\frac{2\delta}{r^{1/2}}+\frac{1}{\delta^2}\leq  
r^{-1}\delta^4+\frac{2}{\delta^2}  \quad\textrm{for }
\delta\in(r^{1/6},r^{1/3}),
\end{eqnarray}
obtained by evaluating $\lambda_m^2+2r^{-1/2}\lambda_m^{-1}$ for  $\lambda_m=\delta^{1/2}$ and $\lambda_m= r^{-1/2}\delta^2$ respectively. To show that $(\ref{1407211722})$ holds it is convenient to operate the change of variable $(r^{1/4},r^{1/2})\ni y:=\delta^{3/2}$ and thus writing $(\ref{1407211722})$ as follows
\begin{eqnarray}\label{1407211726}
y^2(r^{1/2}-2)+2y-r^{1/2}\geq 0\quad\textrm{for } y\in(r^{1/4},r^{1/2})
\end{eqnarray}
which can be easily shown to be true $\forall r\geq 1$.
Then, it is immediate to prove $(\ref{1407211724})$.


\end{document}